\documentclass[letterpaper,10pt]{IEEEtran}
\usepackage{cite}
\usepackage{textcomp}
\usepackage[utf8]{inputenc}
\usepackage[T1]{fontenc}
\usepackage{graphicx}
\usepackage{balance}
\graphicspath{ {./graphics/} }
\usepackage{grffile}
\usepackage{longtable}
\usepackage{wrapfig}
\usepackage{rotating}
\usepackage[normalem]{ulem}
\usepackage{amsmath}
\usepackage{amssymb}
\usepackage{textcomp}
\usepackage{amssymb}
\usepackage{capt-of}
\usepackage{caption}
\usepackage{subcaption}
\usepackage{hyperref}
\usepackage{listings}
\usepackage{float}
\usepackage{optidef}
\usepackage{subfiles}
\usepackage{subfiles}
\usepackage{array}
\usepackage{enumerate}
\usepackage{empheq}
\usepackage{pgfplots}
\usepackage{pifont}
\usepackage{import}
\newcommand{\setR}{\mathbb{R}}
\newcommand{\setN}{\mathbb{N}}
\newcommand{\setRp}{\setR_{\geq 0}}

\newcommand{\ml}[1]{\mathcal{#1}}
\newcommand{\ball}{\mathcal{B}}
\newcommand{\ones}{\boldsymbol{1}}
\newcommand{\conc}[1]{\boldsymbol{#1}}
\newcommand{\abs}[1]{\left| #1 \right|}
\newcommand{\norm}[1]{\lVert #1 \rVert}
\newcommand{\proj}[1]{\mathbb{P}_{#1}}
\newcommand{\indic}[1]{\mathcal{I}_{#1}}

\newcommand{\tderiv}[2]{\frac{d #1}{d #2}}
\newcommand{\jac}{\boldsymbol{\mathrm{J}}}
\newcommand{\clarke}{\jac^c}
\newcommand{\conserv}{\boldsymbol{\mathcal{J}}}
\DeclareMathOperator{\diag}{diag}

\DeclareMathOperator{\inter}{int}
\DeclareMathOperator*{\argmin}{argmin}

\DeclareMathOperator{\dist}{dist}
\DeclareMathOperator{\ncone}{N}

\DeclareMathOperator{\conv}{conv}

\DeclareMathOperator{\sol}{SOL}
\DeclareMathOperator{\vi}{VI}
\DeclareMathOperator{\fix}{fix}

\newcommand{\agents}{\mathcal{N}}
\newcommand{\cy}{\conc{y}}
\newcommand{\cys}{\cy^{\star}}
\newcommand{\com}{\conc{\omega}}
\newcommand{\polyset}{\ml{N}_p}

\newtheorem{stassumption}{Standing Assumption}

\newtheorem{fact}{Fact}

\newtheorem{example}{Example}

\newcounter{algorithm}
\newenvironment{algorithm}[1][]
{	\refstepcounter{algorithm}
	\begin{minipage}{\linewidth}
		\medskip
		\hrule
		\smallskip
		\textsc{Algorithm \thealgorithm}. #1
		\smallskip
		\hrule 
		\smallskip
	\end{minipage}
}
{
	\smallskip
	\hrule width\linewidth\relax
	\smallskip
}

\newcounter{eslalg}

\newtheorem{theorem}{Theorem}
\newtheorem{definition}{Definition}
\newtheorem{remark}{Remark}
\newtheorem{proposition}{Proposition}

\newtheorem{lemma}{Lemma}

\newtheorem{assumption}{Assumption}

\newcommand{\blue}[1]{{\color{black} #1}}

\def\BibTeX{{\rm B\kern-.05em{\sc i\kern-.025em b}\kern-.08em
		T\kern-.1667em\lower.7ex\hbox{E}\kern-.125emX}}
\begin{document}
	\title{
BIG Hype: Best
Intervention in Games via Distributed Hypergradient Descent
	}
	\author{Panagiotis D. Grontas, Giuseppe Belgioioso, Carlo Cenedese, \\ Marta Fochesato, John Lygeros, Florian D\"orfler
		\thanks{This  work  was supported by NCCR Automation, a National Centre of
Competence in Research, funded by the Swiss National Science Foundation (grant number $180545$)
and by the Onassis Foundation - Scholarship ID: F ZQ 019-1/2020-2021. }
		\thanks{Authors are with Automatic Control Laboratory, Department of Electrical Engineering and Information Technology,
        ETH Z\"urich, Physikstrasse 3 8092 Z\"urich, Switzerland. (e-mail:{\tt \{pgrontas, gbelgioioso, ccenedese, mfochesato, jlygeros, dorfler\}@ethz.ch}).
        }
		}
\maketitle
\begin{abstract}
Hierarchical decision making problems, such as bilevel programs and Stackelberg games,
are attracting increasing interest in both the engineering and machine learning communities.
Yet, existing solution methods lack either convergence guarantees or computational efficiency, due to the absence of smoothness and convexity.
In this work, we bridge this gap by designing a first-order hypergradient-based algorithm 
for Stackelberg games
and mathematically establishing its convergence using tools from nonsmooth analysis.
To evaluate the \textit{hypergradient}, namely, the gradient of the upper-level objective,
we develop an online scheme that simultaneously computes the lower-level equilibrium and its Jacobian.
Crucially, this scheme exploits and preserves the original hierarchical and distributed structure of the problem, which renders it scalable
and privacy-preserving.
We numerically verify the computational efficiency and scalability of our algorithm on a large-scale hierarchical demand-response model.
\end{abstract}


\section{Introduction}
\label{sec:introduction}
	A vast array of engineering, social, economic, and financial systems involve a large number of rational decision makers, or agents, interacting with each other and the environment.
	Typically, the resulting global behaviour of these multi-agent systems is inefficient, owing to the self-interested nature of the agents and the lack of coordination.
	It is therefore important to design intervention mechanisms, such as incentives and pricing, that drive the agents to a desirable equilibrium configuration.
	This inherently hierarchical structure appears in such
	diverse fields as
	wireless sensor networks \cite{stack_sensor_net}, energy
	demand-response \cite{drm_miqp} 
	and peer-to-peer trading \cite{stack_energy_trading, le2022parametrized} in smart grids, 
	traffic routing 
\cite{paccagnan2021optimal}, 
	and investment networks \cite{galeotti2020targeting}.
	Recently, hierarchical decision making has attracted also the interest of the machine learning research community, e.g., in the context of hyperparameter optimization \cite{pedregosa_hyperparameter} and meta-learning \cite{rajeswaran2019meta}.
	
	From a mathematical perspective, all these problems find a common interpretation in the paradigm of \textit{bilevel games}. In turn, bilevel games are instances of \textit{mathematical programs with equilibrium constraints} (MPECs) \cite{luo_pang_ralph_1996}
	which are notoriously difficult to solve, since they are, in general, ill-posed
	and violate most of the standard constraint qualifications. %
	A popular solution approach for MPECs, originally introduced in \cite{scholtes_relaxation},
	considers a relaxed version of the equilibrium constraints that makes the problem amenable to standard nonlinear programming solvers.
	The resulting nonlinear program is then iteratively solved while driving the relaxation parameter to zero, hence
	recovering a candidate solution to the original MPEC.
	This seminal work inspired a number of relaxation schemes, which are comprehensively reviewed and numerically compared in \cite{hoheisel2013theoretical}.
Unfortunately, these methods typically suffer from numerical instability when the relaxation parameter approaches zero.
	Another popular solution approach among practitioners is the so-called \textit{big-M reformulation} \cite{FortunyAmat1981ARA} that replaces the equilibrium constraints with linear constraints and auxiliary binary variables. 
	The resulting program can then be solved using off-the-shelf mixed-integer solvers.
	Although provably correct, the poor scalability of mixed-integer solvers makes this approach applicable only to small problems 
	\cite{jara2018study}.
%
%
These shortcomings motivated the development of several \textit{ad hoc} solution approaches tailored to special classes of bilevel games \cite{local_stack_seeking, parise2021analysis} and specific applications \cite{stack_energy_trading, le2022parametrized, paccagnan2021optimal, galeotti2020targeting}.

	
	A conceptually different solution approach based on first-order methods has 
	recently emerged in the machine learning community in the context of hyper-parameter optimization 
	and since then enjoyed considerable success \cite{begio2000gradient}.
However, deploying these first-order methods on bilevel games is challenging since evaluating the ``gradient'' of the upper-level objective, generally known as \textit{hypergradient}, requires computing both the lower-level equilibrium and, crucially, its Jacobian.
	Two schemes have been proposed to overcome this challenge: \textit{implicit differentiation} and \textit{iterative differentiation}. 
\blue{In implicit differentiation, the hypergradient is computed by numerically solving a system of linear equations derived via the implicit function theorem. In practice, this method requires inverting matrices that scale quadratically with the dimension of the lower-level problem, which is computationally impractical for large-scale problems \cite{parise2021analysis, dafermos_sensitivity}. To relax this burden, the implicit equation can instead be solved up to a specified tolerance using the conjugate gradient method \cite{pedregosa_hyperparameter}
	or a truncation of Neumann's series decomposition \cite{lorraine2020optimizing}.
In iterative differentiation, a fixed-point iteration is used to solve the lower-level problem \cite{maclaurin2015gradient}. Then, the hypergradient is approximated by backpropagating derivatives through the solution trajectory using (forward or reverse mode) automatic differentiation \cite{baydin2018automatic}. In the context of bilevel games, this approach was first employed in \cite{li2020end} to design end-to-end interventions sidestepping the heavy computational requirements of implicit differentiation. Nevertheless, this method encounters its own set of challenges. If solving the lower-level problem requires many iterations, the solution trajectory may become too long to efficiently backpropagate derivatives. Motivated by these shortcomings, several works have been recently proposed that aim to approximate the hypergradient computation \cite{li2022fully, yang2021provably, li2023achieving}.} A theoretical and numerical comparison of various implementations of differentiation schemes is presented in \cite{on_hypergradient_complexity}.
	
	Although computational aspects of hypergradient have been extensively studied, convergence guarantees for hypergradient-based methods are typically elusive due to the pathological lack of smoothness and convexity of the underlying problem.
The convergence analysis is further complicated by the inexact evaluation of the hypergradients, which is required to relax its computational burden, especially for large-scale problems.

	In this paper, we solve these technical challenges and design a scalable hypergradient-based algorithm for bilevel games with convergence guarantees.
	 Our contribution is three-fold:
	\begin{enumerate}[(i)]
		\item 
		We develop the first distributed and provably-convergent algorithm for single-leader multi-follower bilevel games, called \textit{BIG Hype}, namely, \underline{B}est \underline{I}ntervention in \underline{G}ames using \underline{Hype}rgradient Descent.
		BIG Hype comprises two nested loops. In the inner loop, an approximation of the lower-level solution and its Jacobian is computed simultaneously by the agents in a distributed manner.
		In the outer loop, the leader uses these approximations to run hypergradient descent.
		To mitigate the inexactness of the hypergradients, we design appropriate termination criteria 
		for the inner loop, and establish convergence by using novel tools of nonsmooth analysis \cite{stochastic_sub_tame}, \cite{nonsmooth_implicit}.
		\item We illustrate how BIG Hype preserves and exploits the original hierarchical and distributed structure of the problem, therefore significantly
		reducing its computational footprint and preserving the privacy of the agents.
		\item We tailor BIG Hype to special subclasses of bilevel games that routinely arise in practice, and derive both sharper theoretical results
		and improved computational performance.
		Specifically, we consider the class of linear-quadratic games with parametric polyhedral and affine constraints, respectively.
		For the latter, we further prove that our algorithm boils down to a single-loop scheme whose computational cost is considerably reduced.
	\end{enumerate}
Finally, we corroborate our theoretical findings via numerical simulations on a hierarchical demand-response model \cite{drm_miqp}. 
		We showcase the scalability of our algorithm by effectively deploying it on huge-scale problem instances, and we compare its performances with state-of-the-art mixed-integer approaches.

	\subsection{Notation and Preliminaries}
	\label{subsec:NP}
	\subsubsection*{Basic Notation}
        We adopt the same basic notation as in \cite{drm_miqp}.
        Given two matrices \( A, B \) of compatible dimensions,
        their horizontal and vertical concatenation is denoted \( [A, B] \) and \( [A;\, B] \), respectively.
	We use \( \norm{\cdot} \) to denote the Euclidean and spectral norm for vector and matrix inputs, respectively.
	The ball with center \( x \in \setR^n \) and radius \( r \) is 
	denoted \( \ball(x;\, r) \).
	For a vector-valued differentiable mapping \( F : \setR^n \times \setR^m \to \setR^p \)
	we denote \( \jac_1 F(x,y) \in \setR^{p \times n} \) and 
	\( \jac_2 F(x,y) \in \setR^{p \times m} \)
	the partial Jacobians of \( F \) with respect to the first and second argument, respectively;
	if \( p = 1 \), we use \( \nabla_1 F, \nabla_2 F \) to denote the partial gradients.
	A subset of \( \setR^n \) is a \textit{semialgebraic} set if it can be expressed
	as a finite Boolean combination of sets of the form
	\( \{ x \in \setR^n \, | \, f(x) > 0 \} \) and
	\( \{ x \in \setR^n \, | \, g(x) = 0 \}  \) where \( f, g \) are polynomial
	functions with real coefficients.
	A function is called semialgebraic if its graph is a semialgebraic set.
    We say that a sequence \( \{ \alpha^k \}_{k \in \setN} \subseteq \setR \) is 
    \textit{nonsummable} if \( \sum_{k=0}^{\infty} \alpha^k = \infty \), and
    \textit{square-summable} if \( \sum_{k=0}^{\infty} (\alpha^k)^2 < \infty \).

	\subsubsection*{Convex Analysis and Operator Theory} 
        We adopt standard operator theoretic notation and definitions from \cite{bauschke2017}.
	Given a mapping \( F : \setR^n \to \setR^n \) and a set \( \ml{Y} \subseteq \setR^n \),
	a solution of the variational inequality problem \( \vi(F, \ml{Y}) \) is a vector \( y^{\star} \in \ml{Y} \)
	such that \( F(y^{\star})^{\top} (y - y^{\star}) \geq 0 \) for all \( y \in \ml{Y} \); the solution
	set of \( \vi(F, \ml{Y}) \) is denoted by \( \sol(F, \ml{Y}) \). The mapping $ F$ is $\mu$-strongly monotone, if $(F(x)\!-\!F(y))^\top (x\!-\!y) \geq \mu \left\| x\!-\!y \right\|^2$ for all $x, y \!\in\! \mathbb R^n$; and $ F$ is $L$-Lipschitz continuous, if $\|F(x)\!-\!F(y)\| \leq L \left\| x\!-\!y \right\|$ for all $x, y \!\in\! \mathbb R^n$.
	For a closed and convex set \( \ml{Y} \subseteq \setR^n \), we denote the projection
	of \( x \in \setR^n \)
	onto \( \ml{Y} \) as \( \proj{\ml{Y}}[x] := \argmin_{y \in \ml{Y}} \norm{y - x} \), 
	and the distance of \( x \) from \( \ml{Y} \) by \( \dist(x; \, \ml{Y}) := \min_{y \in \ml{Y}} \norm{y - x} \).
	The convex hull of \( \ml{Y} \) is denoted \( \conv (\ml{Y}) \); and its normal cone by \( \ncone_{\ml{Y}} : \setR^n \rightrightarrows \setR^n \) which is defined as \( \ncone_{\ml{Y}}(x) = \{ v \in \setR^n \, | \, \sup_{y \in \ml{Y}} v^{\top} (y - x) \leq 0  \} \) if \( x \in \ml{Y} \), and \( \varnothing \) otherwise. 
	
	\subsubsection*{Nonsmooth Analysis}
	Any locally Lipschitz mapping \( F : \setR^n \to \setR^m \) is almost everywhere 
	differentiable (in the sense of Lebesgue measure) by Rademacher's theorem.
	Its \textit{Clarke Jacobian} is a set-valued mapping \( \clarke F : \setR^n \rightrightarrows \setR^{m \times n} \), defined as 
	\(	\clarke F(x) = \conv \big\{ \lim_{k \to \infty} \jac F(x^k) \, | \, x^k \to x, \, x^k \in \Omega \big\}, \)
	where \( \Omega \subseteq \setR^n \) is the full-measure set of points where \( F \) is differentiable \cite{nonsmooth_analysis}.
	Clarke Jacobians generalize the familiar subdifferential from convex analysis to the nonconvex regime.
	
	A recently-introduced generalization of the Clarke Jacobian is the conservative Jacobian \cite{conservative_definition}.
	A set-valued mapping \( \conserv F : \setR^n \rightrightarrows \setR^{m \times n} \) is a \textit{conservative Jacobian}
	if it has a closed graph, is locally bounded, and is nonempty with
	\begin{equation*}
		\tderiv{}{t} F(\rho(t)) \in \conserv F (\rho(t)) \dot{\rho}(t) ~ \text{almost everywhere,}
	\end{equation*}
	where \( \rho : \setR \to \setR^n \) is an absolutely continuous function.
	A function \( F \) that admits a conservative Jacobian is called \textit{path differentiable}.
	Unlike Clarke Jacobians, conservative Jacobians allow for operational calculus, e.g.,
	a nonsmooth implicit function theorem \cite{nonsmooth_implicit}.
	Further, the Clarke and conservative Jacobians satisfy the inclusion
	\( \clarke F(x) \subset \conserv F(x) \).
	
\section{Problem Formulation}
\label{sec:problem_formulation}
We consider hierarchical games consisting of two levels.
	On the lower level, a set of self-interested agents are playing a parametric game.
	On the upper level, a single agent controls the lower-level parameter with the aim of driving the game towards a desirable equilibrium.
	This setting corresponds to a single-leader multi-follower Stackelberg game \cite{von1952theory}.
	
	\subsection{Lower-Level Game}
	At the lower-level game, we consider a group of \( N \) agents, also referred to as followers, and indexed
	by \( i \in \agents := \{1, \ldots, N\} \).
	Each agent \( i \in \agents \) chooses its decision variable \( y_i \) from
	the local feasible set \( \ml{Y}_i(x) \subseteq \setR^{n_i} \) that may be affected
	by the exogenous parameter \( x \in  \setR^m \), namely, the upper-level variable.
	Let \( \cy := (y_i)_{i \in \agents} \in \ml{Y}(x) \) denote the collective strategy
	profile of all agents where
	\( \ml{Y}(x) := \prod_{i \in \agents} \ml{Y}_i(x) \subseteq \setR^n \) and
	 \( n := \sum_{i \in \agents} n_i \).
	The goal of follower \( i \) is to minimize its local objective function
	\( f_i(x, y_i, \cy_{-i}) \) that depends on the parameter \( x \), 
	its own decision \( y_i \), 
	as well as the decisions of the other agents 
	\( \cy_{-i} := (y_j)_{j \in \agents \backslash \{i\}} \in \setR^{n - n_i}  \).
	Overall, the lower-level game is described by the following collection of inter-dependent parametrized optimization
	problems:
	\begin{equation} \label{eq:lower_level}
		(\forall i \in \agents) ~~
		\quad
		\begin{alignedat}{2}
			&\underset{y_i \in \mathbb R^{n_i}}{\mathclap{\mathrm{minimize}}} 
			~~ && ~ f_i(x, y_i, \conc{y}_{-i})
			~\text{ s.t. }~ y_i \in \mathcal{Y}_i(x).\\
		\end{alignedat}
	\end{equation}
	A relevant solution concept for \eqref{eq:lower_level} is the Nash
	equilibrium, that intuitively corresponds to a strategy profile
	in which no agent can reduce its cost by unilaterally deviating from it.
	\begin{definition}
		For any fixed \( x \in \ml{X} \), a strategy profile \( \cys \in \ml{Y}(x) \) is a Nash equilibrium (NE)
		of \eqref{eq:lower_level} if for
		all \( i \in \agents \):
		\begin{equation*}
			f_i(x, y_i^{\star}, \cys_{-i}) \leq f_i(x, y_i, \cys_{-i}), ~ \forall y_i \in \ml{Y}_i(x). \tag*{\( \square \)}
		\end{equation*}
	\end{definition}

\smallskip
	The next \blue{standing assumption\footnote{\blue{Standing assumptions are always valid throughout the paper.}}} on the cost and constraints of the local optimization problems is standard.
	\begin{stassumption} \label{ass:followers_generic_obj_feas}
		For all \( i \in \agents \), and any fixed \( x \in \ml{X} \), \( \cy_{-i} \in \setR^{n - n_i} \),
		the function \( f_i(x, \cdot, \cy_{-i}) \) is convex and continuously differentiable;
		the set \( \ml{Y}_i(x) \) is of the form
		\begin{equation*}
			\ml{Y}_i(x) := \{ y_i \in \setR^{n_i} \, | \, A_i y_i \leq b_i + G_i x, ~ C_i y_i = d_i + H_i x \},
		\end{equation*}
		{\color{black} where
		\( A_i \in \setR^{p_i \times n_i} \) and  \( C_i \in \setR^{r_i \times n_i} \) 
		are full row rank, and
		\( G_i \in \setR^{p_i \times m}, b_i \in \setR^{p_i}, 
		   H_i \in \setR^{r_i \times m}, d_i \in \setR^{r_i} \).}
		\( \ml{Y}(x) \) is nonempty and satisfies Slater's
		constraint qualification.
		\hfill \( \square \)
	\end{stassumption}
	
	Under this assumption, a strategy profile \( \cys \) is a NE of \eqref{eq:lower_level}
	if and only if it solves the parametrized variational inequality \( \vi(F(x,\cdot), \ml{Y}(x)) \) \cite[Prop.\ 1.4.2]{facchinei2003finite},
	where \( F(x, \cdot) \) is the so-called \textit{pseudo-gradient} (PG) mapping, defined as
	\begin{equation} \label{eq:pg_definition}
		F(x, \cy) := (\nabla_{y_i} f_i(x, \cy))_{i \in \agents}.
	\end{equation}

	The following assumption ensures
	existence and uniqueness of a NE for all possible parameters \cite[Th.\ 2.3.3(b)]{facchinei2003finite}.
	\begin{stassumption}\label{ass:pg_stronglym_lipschitz}
		For any fixed \( x \in \ml{X} \), \( F(x, \cdot) \) is \( \mu \)-strongly monotone
		and \( L_F \)-Lipschitz continuous. \hfill \( \square \)
	\end{stassumption}

The unique NE is characterized by the \textit{parameter-to-NE (or solution) mapping} 
	\( \cys : x \mapsto \sol(F(x, \cdot), \ml{Y}(x)) \).
	Next, we postulate some regularity conditions on the PG.
	\begin{stassumption} \label{ass:followers_pg_general}
		The mapping \( F \) is definable\footnotemark, continuously differentiable,
            \footnotetext{\textit{Definable} functions form a vast class of functions that encompasses virtually
    	every objective function employed in optimization or machine learning, e.g.,
    	semialgebraic functions, 
    	exponentials, and logarithms. 
    	See \cite[Def.\ 1.5]{coste2000introduction} for a formal definition and \cite[App.\ A.2]{nonsmooth_implicit}
    	for a practitioner's perspective. 
            If the \( f_i \)'s are definable, so is \( F \), while the inverse is not true \cite[Remark 8]{conservative_definition}.
            }%
		and there exist constants \( L_{JF1}, L_{JF2} \) such that the partial Jacobians
		of \( F \) satisfy
		\begin{align}
			\label{eq:lipshitz_pg1}
			\norm{\jac_1 F(x, \cy) - \jac_1 F(x, \cy')} & \leq L_{JF1} \norm{\cy - \cy'}, \\
			\label{eq:lipshitz_pg2}
			\norm{\jac_2 F(x, \cy) - \jac_2 F(x, \cy')} & \leq L_{JF2} \norm{\cy - \cy'},
		\end{align}
		for any \( \cy, \cy' \in \ml{Y}(x) \).
		\hfill \( \square \)
	\end{stassumption}
	
	A concrete example satisfying this assumption is the class of games in which the cost functions are of the form
	\begin{equation*} \label{eq:simple_general_pg}
		f_i(x, \cy) = \frac{1}{2} y_i^{\top} Q_i(x) y_i + 
		\Big( \textstyle \sum_{j \in \agents} E_{i,j}(x) y_j + E_{i,0} x + e_i \Big)^{\top} y_i,
	\end{equation*}
	where \( Q_i(x) = Q_{i,0} + \sum_{j=1}^{m} Q_{i,j} [x]_j \), with \( [x]_j \)
	denoting the \( j \)-th component of \( x \in \setR^m \) and
	\( Q_{i,j} \in \setR^{n_i \times n_i} \), and
	\( E_{i,j}(x) \) is defined similarly.
	This class is broad and encompasses both linear-quadratic network \cite{parise2021analysis}
	and aggregative games \cite{local_stack_seeking}.

 
	\subsection{Upper-Level Optimization}
	The leader controls the game parameter \( x \), and seeks to solve the following optimization problem:
	\begin{subequations} \label{eq:upper_level_expanded}
		\begin{alignat}{2}
			&\underset{\displaystyle x, \cy}{\mathclap{\mathrm{minimize}}} 
			\quad~ && \quad \varphi(x, \cy) \label{eq:L-OBJF}\\
			& \overset{\hphantom{\displaystyle x , \cy}}{\mathclap{\mathrm{subject~to}}} \quad~
 			&& \quad x \in \ml{X} \\
			& \quad~ && \quad  \cy \in \sol(F(x, \cdot), \ml{Y}(x)), \label{eq:parametric_vi_constraint}
		\end{alignat}
	\end{subequations}
	where \( \ml{X} \subseteq \setR^m \)  is the leader's strategy set
	and \( \varphi \) is the cost function.
	Crucially, \( \varphi \) is affected by \( x \) both explicitly and implicitly
	via the parametric VI constraint \eqref{eq:parametric_vi_constraint}.
	Finally, we postulate regularity of the leader's cost and strategy set.
	\begin{stassumption} \label{ass:leaders_generic}
		The function \( \varphi (x,\cy) \) is definable and continuously differentiable;
		the partial gradients \( \nabla_1 \varphi, \nabla_2 \varphi \) are
		\( L_{\varphi 1}\text{-}, L_{\varphi 2}\text{-}\)Lipschitz continuous, respectively.
		The feasible set \( \ml{X} \) is nonempty, convex, and compact.
		\hfill \( \square \)
	\end{stassumption}

	
	\subsection{Modelling Coupling Constraints} \label{subsec:coupling_constraints}
	Coupling constraints in games as in \eqref{eq:lower_level} are typically used to model applications in which the agents have access to shared resources with limited availability \cite{belgioioso2022distributed}.
	While we do not explicitly consider coupling constraints in \eqref{eq:lower_level}, 
	we can implicitly model them using the parametric feasible set \( \ml{Y}(x) \).
	Intuitively, the leader can locally restrict each agent's access to the shared resource such that the overall resource availability is respected, \blue{as explained in \cite[\S~3.2]{nabetani2011parametrized}. 
Formally, consider the game in \eqref{eq:lower_level} with additional coupling constraints of the form
	\begin{equation} \label{eq:coupling_constraints}
		\textstyle \sum_{i \in \agents} \Theta_i y_i \leq \theta,
	\end{equation}
	where \( \Theta_i \in \setR^{n_{\theta} \times n_i} \) and
	\( \theta \in \setR^{n_{\theta}} \). By endowing the leader with the auxiliary  variables \( \hat{\theta}_i \in \setR^{n_{\theta}} \), for all $i \in \mathcal{N}$, we can equivalently recast the coupling constraints \eqref{eq:coupling_constraints} as
	\begin{subequations}
	    \begin{align}
		\Theta_i y_i & \leq \hat{\theta}_i, \quad \forall i \in \mathcal{N} \label{eq:coupling_per_agent} \\
		\textstyle \sum_{i \in \agents}  \hat{\theta}_i & \leq \theta \label{eq:coupling_leader}.
	\end{align}
	\end{subequations}
	Hence, each inequality in \eqref{eq:coupling_per_agent} becomes a local constraint embedded in the corresponding set \( \ml{Y}_i(x) \),
	whereas \eqref{eq:coupling_leader} can be incorporated in the leader's constraints set \( \ml{X} \). This reformulation gives the leader increased interventional flexibility. For example, social welfare criteria can be prescribed by augmenting the leader's objective function \eqref{eq:L-OBJF} with terms that promote a "fair" distribution of shared resources. In Section~\ref{sec:numerical_simulations}, we showcase a practical example of such reformulation.}

\section{A Distributed Hypergradient Algorithm}
\label{sec:algorithm}
	We approach \eqref{eq:upper_level_expanded} by eliminating the 
	parametric VI constraint \eqref{eq:parametric_vi_constraint} using the single-valued solution mapping \( \cys(\cdot) \),
yielding the nonsmooth nonconvex optimization problem
	\begin{equation} \label{eq:upper_level_substi}
		\begin{alignedat}{2}
			&\underset{\displaystyle x \in \ml{X}}{\mathclap{\mathrm{min}}} 
			\quad~ && \varphi(x, \cys(x)) =: \varphi_e(x). \\
		\end{alignedat}
	\end{equation}

	The core of our algorithm is simple:
	compute the gradient of \( \varphi_e \), commonly referred to as the \textit{hypergradient}, and solve \eqref{eq:upper_level_substi}
	using projected hypergradient descent.
	Whenever \( \cys(\cdot) \) is differentiable at $x$, an expression for the hypergradient can be derived using the chain rule as follows:
	\begin{equation}
	\label{eq:hypergrad}
		\nabla \varphi_e(x) = \nabla_1 \varphi(x, \cys(x)) + \jac \cys(x)^{\top} \nabla_2 \varphi(x, \cys(x)).
	\end{equation}
	A similar expression holds when \( \cys(\cdot) \) is nondifferentiable at $x$, where standard Jacobians
	are replaced with elements of the conservative Jacobian \cite{nonsmooth_implicit}.
        We refer to the proof of \autoref{th:convergence_lqg} in Appendix \ref{proof:convergence_lqg} for a formal analysis of this case.
	Obtaining \( \nabla \varphi_e(x) \) requires the exact evaluation of \( \cys(x) \) as well as its Jacobian \( \jac \cys(x) \), also known as \textit{sensitivity}, which can be computationally prohibitive for large-scale games \cite{parise2021analysis}. Instead, we rely on approximations and enforce accuracy bounds to mitigate the errors arising from inexactness of the estimates of \( \cys(x) \) and \( \jac \cys(x) \). To alleviate the computational burden further and to preserve \blue{data privacy of the followers}, we develop an algorithm that allows evaluating both \( \cys(x) \) and \( \jac \cys(x) \) simultaneously and in a distributed way.
	
	The proposed scheme is summarized in \autoref{alg:general_hypergradient} and consists
	of two nested loops.
	In the inner loop, the agents receive the current leader's strategy \( x \) and update their local estimates of the equilibrium and the sensitivity, until an appropriate termination criterion is met.
In the outer loop, the leader computes the approximate hypergradient \( \widehat{\nabla} \varphi_e \) by gathering the local equilibrium and sensitivity estimates, respectively \( \cy^k \) and \( s^k \), and then performs 
	a projected hypergradient step. \blue{The discussion on the role of $\alpha^k$ is delayed to \autoref{subsec:outer_loop_an}.} 

    First, we introduce and motivate the basic iterations that will form the core of the inner loop scheme; the formal algorithm is summarized in \autoref{alg:innner_general_setup} and analysed in \autoref{sec:convergence_analysis}.
    \blue{In \autoref{alg:innner_general_setup}, we will denote the inner loop iterates as \( \tilde{\cy} \) and \( \tilde{s} \), whereas
    \( \bar{\cy} \) and \( \bar{s} \) will indicate the respective output.}
	\begin{figure} [t]
		\begin{algorithm}[BIG Hype]
		\textbf{Parameters}:
			step sizes \( \{ \alpha^k, \beta^k \}_{k \in \setN} \), tolerances \( \{\sigma^k\}_{k \in \setN} \).\\[.2em]
				\textbf{Initialization}: $k\leftarrow 0$, 
				\( x^k \in \ml{X},~ \cy^k \in \setR^n, ~ s^k \in \setR^{m \times n} \).
				\\[.2em]
				\textbf{Iterate until convergence:}\\
$
					\left \lfloor
					\begin{array}{l}
						\text{Leader's projected hypergradient step:} \\[.1em]
						\left|
						\begin{array}{l}
						\widehat{\nabla} \varphi_e^k =
						\nabla_1 \varphi(x^k, \cy^{k}) + (s^k)^{\top} \nabla_2 \varphi(x^k, \cy^{k}) \\
						x^{k+1} = x^k + \beta^k \big( \proj{\mathcal{X}}[x^{k} - \alpha^k \widehat{\nabla} \varphi_e^k] - x^k\big) \\
						\end{array}
						\right.\\[1.15em]
				   %
\text{Followers' Approximate NE and Sensitivity Learning:}\\[.2em]
\left|
				    \begin{array}{l}
				     (\cy^{k+1}, \ s^{k+1}) =	\\
				     ~~~~
 \text{\textbf{Distributed Inner Loop}} \\[.1em]
~~~~			    \left\lfloor
			    \begin{array}{l}
				     \text{Input: }
x^{k+1},\, \cy^k,\, s^k,\, \sigma^k\\
				        \text{Output: }
(\bar \cy, \, \bar s ) \, \text{ such that }\\
{\small
~\max\left\{\|\bar \cy\!-\!\cy^\star(x^{k+1})\|, \, \| \bar s\!-\!\jac \cy^\star(x^{k+1})\| \right\} \leq \sigma^k
}
				    \end{array}
				    \right. \end{array}
				    \right. \vspace{.25em}
				    \\
					k \leftarrow k+1
					\end{array}
					\right.
$
				\label{alg:general_hypergradient}
		\end{algorithm}
	\end{figure}
	
	\subsection{Equilibrium and Sensitivity Learning}

	In the inner loop, whose iterations are indexed by \( \ell \in \setN \), we estimate \( \cys(x) \) using a fixed-point iteration of the form
	\begin{equation} \label{eq:y_update}
		(\forall \ell \in \setN) \quad 
		\tilde{\cy}^{\ell+1} =
h(x, \tilde{\cy}^{\ell}), 
	\end{equation}
	where $h(x,\cdot)$ is the \textit{projected pseudo-gradient} (PPG) map, i.e.,
	\begin{equation}
		h(x, \cdot) := \proj{\ml{Y}(x)}[~\cdot - ~\gamma F(x, \cdot)], \quad \text{with } \gamma > 0.
	\end{equation}
	It can be shown that \( \cys(x) \) is equivalently characterized as a fixed point of the PPG mapping, namely, for any given \( x \in \ml{X} \)
\begin{equation*}
\cy \in \sol(F(x, \cdot), \ml{Y}(x)) \Leftrightarrow \cy = h(x, \cy).
\end{equation*}
	\blue{Additionally, under \autoref{ass:followers_generic_obj_feas}, for any choice of
	\( \gamma < 2 \mu / L_F^2 \), the mapping \( h(x, \cdot) \) is a contraction
	with constant \( \eta := \sqrt{1 - \gamma(2 \mu - \gamma L_F^2)} < 1 \) and, thus, 
	the sequence generated by \eqref{eq:y_update}
	converges to \( \cys(x) \) linearly with rate \( \eta \) \cite[Prop.~26.16]{bauschke2017}.}
        Finally, we note that the Cartesian structure of \( \ml{Y} \) allows for
	an agent-wise decomposition of \eqref{eq:y_update}:
	\begin{equation} \label{eq:y_update_distributed}
		(\forall \ell \in \setN) (\forall i \in \agents) \quad 
		\tilde{y}_i^{\ell+1} = h_i(x, \tilde{\cy}^\ell)
	\end{equation}
	where \( h_i(x, \cy) := \proj{\ml{Y}_i(x)}[y_i - \gamma F_i(x, \cy)] \) and $F_i$ is the $i$-th component of the PG, namely, $F_i(x, \cy) := \nabla_{y_i} f_i(x,\cy)$.

\smallskip
        To estimate \( \jac \cys(x) \), one can use the fixed-point iteration
	\begin{equation} \label{eq:s_update_nominal}
		(\forall \ell \in \setN) \quad \hat{s}^{\ell+1} = \jac_2 h(x, \cys(x)) \hat{s}^\ell
		+ \jac_1 h(x, \cys(x))
	\end{equation}
	which is intuitively motivated by differentiating \eqref{eq:y_update} at its fixed point \( \cys(x) \). 
	%
	%
	In \autoref{sec:convergence_analysis}, we formally prove that \eqref{eq:s_update_nominal} converges to \( \jac \cys(x) \) using the implicit function theorem and the fact that \( h(x, \cdot) \) is contractive. In practice, implementing \eqref{eq:s_update_nominal} is problematic as it requires 
	perfect knowledge of \( \cys(x) \), which can only be obtained via an infinite number of PPG iterations \eqref{eq:y_update}.
	Instead, we propose an online version of \eqref{eq:s_update_nominal} that approximates
	\( \cys(x) \) using the latest iterate of \eqref{eq:y_update}, i.e.,
	\begin{equation} \label{eq:s_update_perturbed}
(\forall \ell \in \setN)	\quad \tilde{s}^{\ell+1} = \jac_2 h(x, \tilde{\cy}^{\ell+1}) \tilde{s}^\ell
		+ \jac_1 h(x, \tilde{\cy}^{\ell+1}).
	\end{equation}
	This sensitivity learning scheme is not a standard fixed-point iteration, 
	due to the update rule depending on \( \tilde{\cy}^{\ell+1} \), 
	but can be seen as a perturbed version of the nominal iteration \eqref{eq:s_update_nominal}.
	
	Remarkably, \eqref{eq:s_update_perturbed} preserves the original distributed structure of the problem, thus making it amenable to a distributed implementation.
	This distributability is evident by decomposing the Jacobian estimate \( \tilde{s} \) in \( N \) row blocks \( \tilde{s}_i \in \setR^{n_i \times m} \),
	each of which is owned by the corresponding agent \( i \in \agents \).
	Then, we can express the sensitivity update \eqref{eq:s_update_perturbed} agent-wise as 
	\begin{equation} \label{eq:s_update_distributed}
		(\forall i \in \agents) ~~~ 
		\tilde{s}_{i}^{\ell+1} = \jac_2 h_i(x, \tilde{\cy}^{\ell+1}) \, \tilde{s}^\ell
		+ \jac_1 h_i(x, \tilde{\cy}^{\ell+1}).
	\end{equation}
	\blue{
By defining \( g_i(x, y) := \proj{\ml{Y}_i(x)}[y] \), for ease of exposition, we can write the partial Jacobians in \eqref{eq:s_update_distributed} as\footnote{The explicit dependence of $\mathbf{J}_1 h_i$ and $\mathbf{J}_2h_i$ on $x,\cy$ is omitted for brevity.}%
	\begin{subequations}%
	\label{eq:ppg_jac}
    	\begin{align}
    		\jac_1 h_i & = \jac_1 g_i(x, y_i - \gamma F_i(x, \cy)) \label{eq:ppg_jac_1_i}  \\
    		&  ~~~ - \gamma \jac_2 g_i(x, y_i - \gamma F_i(x, \cy)) \jac_1 F_i(x, \cy), \notag \\
    		\jac_2 h_i & = \jac_2 g_i(x, y_i - \gamma F_i(x, \cy)) (I - \gamma \jac_2 F_i(x, \cy)),
    		\label{eq:ppg_jac_2_i}
    	\end{align}
	\end{subequations}	
where both expressions require differentiating through \( g_i(x, y) \) which is a projection onto a parametric polyhedron. Analytical properties and computational methods for this operation have been studied in \cite{optnet_arxiv}, and are included in Appendix~\ref{app:diff_proj} for the sake of completeness. In practice, computing the partial Jacobians $\jac_1 g_i(x,y)$ and $\jac_2 g_i(x,y)$ require solving a linear system of equation, as explained in Appendix~\ref{app:diff_proj}.
	}	
	
	We stress that both the equilibrium \eqref{eq:y_update_distributed} and sensitivity \eqref{eq:s_update_distributed} updates can be computed by the agents using solely locally-accessible information\footnote{In partial-decision information settings \cite{belgioioso2022distributed}, namely, when the agents do not have access to all the decisions of the agents that directly influence their cost functions, the updates \eqref{eq:y_update_distributed} and \eqref{eq:s_update_distributed} can be augmented with auxiliary variables and consensus dynamics to retain convergence while respecting the information structure of the game. We leave this extension as future work.}, namely, the decisions of the leader and the other agents that directly influence their cost function. 	
	
	\blue{
\begin{remark}
When \(h_i(x,\cdot)\) lacks differentiability at \({\cy}\), \eqref{eq:ppg_jac_1_i} and \eqref{eq:ppg_jac_2_i} produce matrices $S_{1,i}$ and $S_{2,i}\), respectively, which are elements of the Clarke Jacobian \(\clarke h_i(x, {\cy})\), rather than the standard Jacobians \(\jac_1 h_i(x, {\cy}), \jac_2 h_i(x, {\cy})\). We postpone a formal discussion of this case to \autoref{sec:convergence_analysis}.
{\hfill $\square$}
\end{remark}	
	}

	%
\begin{figure}[t]
		\begin{algorithm}[Distributed Inner Loop]%
\label{alg:innner_general_setup}%
\textbf{Parameters}: step size \( \gamma \),	contraction constant \( \eta \). \\
\textbf{Input}:
			\(  x, \cy, s , \sigma\).
			\\
\textbf{Initialization:}  \(\ell \leftarrow 0, \ \textrm{termination}=\textrm{false},   \)\\[.2em]
\( ~~~ \tilde s^\ell = s, \ \,  \tilde \cy^\ell =\begin{array}{l}
						\left\{
						\begin{array}{l l}
\textbf{Warmstart}\,(x, \cy, \sigma),
& \text{  General}\\
\cy, & \text{  LQG, LQSG} 				
						\end{array}
						\right.
\end{array}
\)
\\[.5em]
\textbf{Iterate until termination} \\[0.5em]
$
					\begin{array}{l}
						\left\lfloor
						\begin{array}{l}
							\text{For all agents \( i \in \agents \):} \\ 
							\left\lfloor
							\begin{array}{l}
							\text{Equilibrium seeking step:}\\
							\left|
							\begin{array}{l}
		\tilde{y}_i^{\ell+1} = h_i(x,\tilde{\cy}^{\ell})
							\end{array}
									\right.
\\[.5em]
                            \text{Jacobian update:} \\
                         \left| 
 \begin{array}{l l}
                                \begin{array}{l l}
                                S_{1,i} \in \clarke_1 h_i(x, \tilde{\conc{y}}^{\ell+1}), & \text{ 
 via \eqref{eq:ppg_jac_1_i}}\\
                           S_{2,i} \in \clarke_2 h_i(x, \tilde{\conc{y}}^{\ell+1}), & \text{  via \eqref{eq:ppg_jac_2_i}} 
                                \end{array}
                            \end{array}
                             \right.
\\[.5em]
		\text{Sensitivity learning step:}\\
							\left|
							\begin{array}{l}
			\tilde{s}_i^{\ell+1} = 
			S_{2, i} \tilde{s}^\ell + S_{1, i}	
							\end{array}
									\right.	
								\vspace*{.5em}
\end{array}
\right. \\[1.5em]
\textrm{termination} = \\
\; \left\{
\begin{array}{l l}
\max\big\{(\eta)^\ell,  \sum_{j=0}^\ell (\eta)^{\ell-j} \|\tilde{\cy}^{j+1} \!-\! \tilde{\cy}^{j} \| \big\} \leq \sigma, & \text{General} \\
\max \big\{\|\tilde{\cy}^{\ell+1}\! -\! \tilde{\cy}^{\ell} \|, \|\tilde{s}^{\ell+1}\! -\! \tilde{s}^{\ell} \| \big\} \leq \sigma, & \text{LQG}\\
\textrm{true}, & \text{LQSG} \\
\end{array}
\right.\\[2em]
							\ell \leftarrow \ell + 1  \\
						\end{array}
						\right. 
					\end{array}	
$\\[.2em]
\textbf{Output}: \( \bar \cy= \tilde{\cy}^{\ell}, \bar s= \tilde{s}^{\ell} \).
		\end{algorithm}
	\end{figure}	
	\begin{figure}[t]
		\begin{algorithm}[Warmstart Loop]%
\label{alg:warmstart}%
\textbf{Input}:
			\(  x, \ \cy,\ \sigma \).
			\\
\textbf{Initialization:}  \(j \leftarrow 0, \ \textrm{termination}=\textrm{false}, \ \tilde \cy^j =
\cy. \)\\
\textbf{Iterate until termination} \\[.5em]
$
					\begin{array}{l}
						\left\lfloor
						\begin{array}{l}
							\text{For all agents }  i \in \agents: \\ 
\left\lfloor
\begin{array}{l}
\tilde{y}_i^{j+1} = h_i(x,\tilde{\cy}^{j})
\end{array}
\right. \\[.5em]
\textrm{termination} = 
\left(\|\tilde{\cy}^{j+1}\! -\! \tilde{\cy}^{j} \| \leq \sigma\right)\\[.5em]
							j \leftarrow j + 1  \\
						\end{array}
						\right. 
					\end{array}	
$\\[.2em]
\textbf{Output}: \( \bar \cy= \tilde{\cy}^{j}\).
		\end{algorithm}
	\end{figure}	
	\subsection{Preview of the Convergence Results}
	In our convergence analysis, we consider the general problem setting
	prescribed by Standing Assumptions \ref{ass:followers_generic_obj_feas}-\ref{ass:leaders_generic},
	as well as two popular problem classes that arise by imposing additional structure on the lower-level game.
	The first problem class is LQGs, that satisfies the following additional assumption.
	\begin{assumption} \label{ass:affine_pg}
		For each agent \( i \in \agents \),
		the cost function is
		\begin{equation}
					f_i(x,\cy) = \frac{1}{2} y_i^{\top} Q_i y_i + 
			\big( \textstyle \sum_{j \in \agents} E_{i,j} y_j + E_{i,0} x + e_i \big)^{\top} y_i,
		\end{equation}
		where \( Q_i \in \mathbb S_{\succ 0}^{n_i } \), \( E_{i,j} \in \setR^{n_i \times n_j}, E_{i,0} \in \setR^{n_i \times m},  e_i \in \setR^{n_i} \). 
		\hfill \( \square \)
	\end{assumption}
 
	
The second problem class considered is a subset of LQGs where the feasible action set of each agent is affine.
	\blue{This class is termed \textit{linear-quadratic subspace games} (LQSG) and satisfies Assumption~\ref{ass:affine_pg} as well as the following assumption.}
	\begin{assumption} \label{ass:affine_subspace_con}
		For all \( i \in \agents \) and \( x \in \ml{X} \), the local
		feasible set is given by
		\( \ml{Y}_i(x) := \{ y_i \in \setR^{n_i} ~|~ C_i y_i = d_i + H_i x \} \).
		\hfill \( \square \)
	\end{assumption}
 
Notably, LQSGs admit an affine solution mapping of the form \( \cys(x) = W x + w \), for some \( W \in \setR^{n \times m}\) and \(w\in \setR^{n} \).	

	We tailor the inner loop of the general scheme presented in \autoref{alg:general_hypergradient}
	 to exploit the structure of each problem class.
	The three resulting algorithms are jointly summarized in Alg.~\ref{alg:innner_general_setup}, and enjoy different theoretical and computational characteristics that
	are briefly outlined in \autoref{table:convergence_results}.
	\begin{table*}[htbp]
	\centering
	\renewcommand{\arraystretch}{1}
	\begin{tabular}{c||c|c|c|c|c|c|c}
		Problem Setup   & Assumptions & Convergence & \blue{ \( \{\alpha_k\}_{k \in \setN} \) }&  \blue{\( \{\beta_k\}_{k \in \setN} \) }& Loop & Warmstart-free & Inner Termination   \\ \hline
		General	 & A3 & Theorem 1 & Vanishing & Constant & Double  & \ding{55}  &
		\textit{A priori} error bound   \\ \hline
		LQG & A1 \& A3
		&  Theorem 2 & Vanishing & Constant & Double &  \ding{51} &
		\textit{A posteriori} error bound   \\ \hline
		LQSG 
		+ convex \( \varphi \)	& A1 \& A2
		& Theorem 3 & Constant & Vanishing
		& Single & \ding{51} &  Perform one iteration \\ \hline
		\renewcommand{\arraystretch}{1}
		LQSG 
		+ strongly convex \( \varphi \) & A1 \& A2
		& Theorem 4 & Constant & Constant & Single & \ding{51} & Perform one iteration \\
	\end{tabular} 
	\caption{Summary of the different problem settings, corresponding inner loop algorithms and their properties.}
	\label{table:convergence_results}
\end{table*}
	
	Notably, the computational properties of our algorithm improve the more assumptions
	are imposed on the problem.
	Specifically, our most general problem setup requires the inner loop to be preceded by a warmstarting phase, summarized in Alg.~\ref{alg:warmstart}, in which the NE is estimated while the sensitivity is not updated. Furthermore, the inner loop termination requires evaluating an \textit{a priori} error bound 
	which is rather conservative.
	By contrast, LQGs utilize an \textit{a posteriori} error bound which is more accurate
	and, hence, typically results in fewer inner iterations.
	For LQSGs with convex \( \varphi \) we can further show that one single inner iteration is sufficient to ensure
	convergence. This renders our algorithm a single-loop scheme and significantly
	decreases its computational burden.
	Finally, a constant step size can be employed for LQSGs with strongly convex \( \varphi \) which
	allows for faster convergence rates.

	\subsection{Scalability in Aggregative Games}
	A practical concern regarding \autoref{alg:general_hypergradient} is the scalability of computational
	and memory requirements with respect to the number of followers.
	In this section, we discuss a common and widely studied game setup,
	corresponding to \textit{aggregative games} \cite{belgioioso2022distributed},
	for which the computational cost of BIG Hype is independent of the number of followers.
	In aggregative games, each
	agent is affected by some aggregate quantity that depends on all the agents' decisions. {\color{black} Based on the dependence on the aggregate function, aggregative games live across the different problem classes. In particular, LQG and LQSG games are
    aggregative games with an affine form of the aggregation function.}
	Here we consider linear aggregates in which the cost function of each agent has the form \( f_i(x, y_i, \sigma(\cy)) \), where
	\begin{equation}
		\sigma(\cy) := \textstyle \sum_{j \in \agents} K_j y_j \in \setR^{\bar{n}}
	\end{equation}
	is the \textit{aggregate variable} of dimension \( \bar{n} \), and \( K_i \in \setR^{\bar{n} \times n_i} \).
	In this case, the PPG is a mapping of the form \( h_i(x, y_i, \sigma(\cy)) \)
	and the update rule \eqref{eq:s_update_perturbed} becomes
	\begin{align}
		\tilde{s}_i^{\ell+1} = 
		& \underbrace{\jac_2 h_i(x, y_i^\ell, \sigma(\cy^\ell))}_{ \in \setR^{n_i \times n_i}} \tilde{s}_i^\ell
		+ \underbrace{\jac_3 h_i(x, y_i^\ell, \sigma(\cy^\ell))}_{ \in \setR^{n_i \times \bar{n}}} \sum_{j \in \agents} K_j    
            \tilde{s}_j^\ell
		\notag \\
		&~~~~~+\jac_1 h(x, y_i^\ell, \sigma(\cy^\ell)), \label{eq:s_update_aggregative}
	\end{align}
        \blue{where the first, second, and third term describe how the agent is influenced by its own decision,
        the aggregate quantity, and the leader's decision, respectively.}
	To implement \eqref{eq:s_update_aggregative} in a scalable manner, 
        we assume that the leader gathers \( K_i \tilde{s}_i^\ell \) from
	all agents and then broadcasts the aggregate sensitivity \( \sum_{i \in \agents} K_i \tilde{s}_i^\ell \).
	Under this communication pattern, the computational cost of \eqref{eq:s_update_aggregative} is independent of 
	\( N \).
	Additionally, if the leader's cost function admits an aggregative structure of the form 
	\( \varphi(x, \sigma(\cy)) \) -- as is often the case in practice \cite{drm_miqp}, \cite{local_stack_seeking} --
	then the hypergradient is given by
	\begin{align*}
		\nabla \varphi_e(x) = & \nabla_1 \varphi(x, \sigma(\cys(x))) \\
		& ~~~~ + \big( \textstyle \sum_{i \in \agents} K_i s_i \big)^{\top} \nabla_2 \varphi(x, \sigma(\cys(x))),
	\end{align*}
	and possesses the same scalability properties.
	

\section{Convergence Analysis}
\label{sec:convergence_analysis}
   In this section, we elaborate on specific algorithms for the inner loop of \autoref{alg:general_hypergradient}, 
    as well as for the two special problem classes, and we establish convergence guarantees. 
    A prominent part of our analysis concerns the sensitivity learning scheme \eqref{eq:s_update_perturbed}, 
    which is tightly bound to the smoothness analysis of parametric polyhedral projection operators.
    The proofs of all technical statements are provided in the appendix.

	\subsection{Inner Loop Analysis}
	\label{subsec:inner_loop_an}
	 \subsubsection{General Games}
%
    We begin our analysis by investigating the differentiability of the parametric projection mapping
	\begin{equation} \label{eq:projection_as_qp}
	\proj{\ml{Y}(x)}[\cy] = 
			\underset{ z \in \ml{Y}(x)}{\mathclap{\mathrm{argmin}}} 
					\quad~ \frac{1}{2} z^{\top} z - \cy^{\top} z 
		\end{equation}
	which constitutes a multiparametric quadratic program and admits a 
	piecewise affine (PWA) solution with respect to \( \cy \) \cite[Th.~1]{mpqp_results}.
	In other words, there exists a \textit{polyhedral partition}
	of \( \setR^n \), namely, a finite family \( \{ \tilde{\ml{Y}}_i \}_{i \in \polyset} \) of polyhedral sets\footnote{The partitions satisfy \( \bigcup_{i \in \ml{N}_p}^{} \tilde{\ml{Y}}_i = \setR^n \) and 
	\( \inter \tilde{\ml{Y}}_i \bigcap \tilde{\ml{Y}}_j = \varnothing \)
	for \( i \neq j \in \polyset \), and their dependence on \( x \) is suppressed for ease of notation.} \( \tilde{\ml{Y}}_i \),
	indexed by \( i \in \polyset := \{1, \ldots, N_p\} \), such that
	\( \proj{\ml{Y}(x)}[\cdot] \) is an affine function on each partition element \( \tilde{\ml{Y}}_i \) and,
	therefore, is continuously differentiable on its interior, \( \inter \tilde{\ml{Y}}_i \).
	Differentiability of \( \proj{\ml{Y}(x)}[\cdot] \) fails only on the boundary of \( \tilde{\ml{Y}}_i \) which, however,  is a lower-dimensional subspace. 
	
	More formally, let \( \com(x, \cy) := \cy - \gamma F(x, \cy) \) be the PG step mapping, and 
	\( \ml{P}(x) := \{ i \in \polyset \, | \, \com(x, \cys(x)) \in \tilde{\ml{Y}}_i \} \) be the set of indices of the active partition elements at \(\com(x, \cys(x)) \).
        The next fact follows directly from the previous paragraph.
        \begin{fact} \label{fact:diff_equivalence}        For any fixed \( x \in \ml{X} \), the following are equivalent:%
            \begin{enumerate}[(i)]%
            \item \( \ml{P}(x) \) is a singleton;
            \item \( \com(x, \cys(x)) \in \inter \tilde{\ml{Y}}_i \) for some unique \( i \in \ml{N}_p \);
            \item \( \proj{\ml{Y}(x)}[\cdot] \) is continuously differentiable at \( \com(x, \cys(x)) \).
            \hfill \( \square \)
            \end{enumerate}
        \end{fact}

	With these observations in place, we derive a condition 	for differentiability of \( \cys(\cdot) \) and establish convergence of the nominal sensitivity learning
	iteration \eqref{eq:s_update_nominal}.

	\begin{proposition} \label{prop:s_nominal_convergence}
		Fix \( x \in \ml{X} \) and assume \( \ml{P}(x) \) is a singleton. 
		Then, \( \cys(\cdot) \) is continuously differentiable at \( x \).
		Further, the sequence \( \{\hat{s}^\ell\}_{\ell \in \setN} \) generated by \eqref{eq:s_update_nominal}
		converges to \( \jac \cys(x) \).
		\hfill \( \square \)
	\end{proposition}
	
	Next, we extend our convergence result to the perturbed iteration \eqref{eq:s_update_perturbed}.
	To this end, we first bound the error in our update rule \eqref{eq:s_update_perturbed}
	by the distance between \( \tilde{\cy}^{\ell} \) and \( \cys(x) \).
	\begin{lemma} \label{lemma:s_perturbed_bound}
		Fix \( x \in \ml{X} \) and assume \( \ml{P}(x) \) is a singleton. 
		Then, there exists \( \bar{\ell} \in \setN \) such that
		for all \( \ell \in \setN \), with \( \ell \geq \bar{\ell} \), the following bounds hold:
        \begin{align}
			\text{(i)} ~ &  \norm{\jac_1 h(x, \tilde{\cy}^\ell) - \jac_1 h(x, \cys(x))} 
                < \gamma L_{JF1} \norm{\tilde{\cy}^\ell - \cys(x)} \notag \\
			\hfill \hphantom{. \square}
			\text{(ii)} ~ &  \norm{\jac_2 h(x, \tilde{\cy}^\ell) - \jac_2 h(x, \cys(x))} 
                < \gamma L_{JF2} \norm{\tilde{\cy}^\ell - \cys(x)}. \notag
			\hfill 
        \end{align}
        
        \vspace*{-.3em}
        {\hfill $\square$}
	\end{lemma}
        \medskip
	
	The condition \( \ell \geq \bar{\ell} \) is a natural consequence of
	the discontinuity of \( \jac h(x, \cdot) \) along the boundary of the polyhedral partition,
	due to the discontinuity of \( \jac \proj{\ml{Y}(x)}[\cdot] \).
	Intuitively, to accurately approximate \( \jac h(x, \cys(x)) \) it must hold
	that \( \cys(x) \) and \( \tilde{\cy}^\ell \) lie on the same partition \( \tilde{\ml{Y}}_i \)
	such that the effect of discontinuities is eliminated.
	Since \eqref{eq:y_update} is a contraction, this can be achieved after finitely many inner-loop iterations.
	\begin{proposition} \label{prop:s_general_convergence}
		Fix \( x \in \ml{X} \) and assume \( \ml{P}(x) \) is a singleton. 
		The sequence \( \{\tilde{s}^\ell \}_{\ell \in \setN} \) generated by \eqref{eq:s_update_perturbed}
		converges to \( \jac \cys(x) \).
	
		\hfill \( \square \)
	\end{proposition}

	Next, we derive bounds on the convergence error of \( \tilde{\cy}^\ell \) and \( \tilde{s}^\ell \)
	that will allow us to control the hypergradient approximation error.
	For the equilibrium error, we recall that \eqref{eq:y_update} satisfies the following
	\textit{a posteriori} error bound \cite[Th.\ 1.50(v)]{bauschke2017}:
	\begin{equation} \label{eq:y_bound}
		\norm{\tilde{\cy}^\ell - \cys(x)} \leq \frac{1}{1 - \eta} \norm{\tilde{\cy}^{\ell+1} - \tilde{\cy}^\ell}.
	\end{equation}
	The sensitivity error is quantified in the following lemma.	
	\begin{lemma} \label{lemma:s_bound_general}
		Fix \( x \in \ml{X} \) and assume that \( \ml{P}(x) \) is a singleton.
		If \( \com(x,\tilde{\cy}^0) \in \tilde{\ml{Y}}_i \) for the unique \( i \in \ml{P}(x) \),
		then there exist constants \( B_{ps}, B_{ns} > 0\) such that the iterates
		of \eqref{eq:s_update_perturbed} satisfy:%
		\begin{equation} \label{eq:s_bound_general}
			\norm{\tilde{s}^\ell - \jac \cys(x)} \leq 
			B_{ps} \sum_{j=0}^{\ell-1} \eta^{\ell-1-j} \norm{\tilde{\cy}^{j+1} - \tilde{\cy}^j} +
			B_{ns} \eta^\ell.
		\end{equation}
		\vspace*{-2em}
		
		{\hfill $\square$}
	\end{lemma}
	
	Note that the condition \( \com(x,\tilde{\cy}^0) \in \tilde{\ml{Y}}_i \) holds without loss of generality, since it can be enforced by running a finite number of PPG steps before deploying sensitivity learning.

	Based on the previous analysis, we propose the inner loop algorithm 
	summarized in \autoref{alg:innner_general_setup} (General), that comprises two phases.
	In the initial phase, summarized in Alg.~\ref{alg:warmstart}, we run only PPG updates until the condition of \autoref{lemma:s_bound_general} is met. This is achieved by enforcing the equilibrium error, which is upper bounded by 
	\( \norm{\tilde{\cy}^{\ell+1} - \tilde{\cy}^{\ell}} \), to be sufficiently small such that
	\( {\tilde \cy^{\ell+1}} \) lies in a partition \(\tilde{\ml{Y}}_j \) for \( j \in \ml{P}(x) \).
	In the second phase, PPG and sensitivity learning steps are performed simultaneously.
	Note that the choice of the termination criterion and the error
	bounds \eqref{eq:y_bound} and \eqref{eq:s_bound_general} guarantee that both the equilibrium and sensitivity errors are smaller than the required tolerance.

\smallskip
	\subsubsection{Linear-Quadratic Games (LQGs)}
	\label{subsec:LQGS}
	Here, we focus on the class of LQGs and extend the convergence results to the nonsmooth domain.
	Firstly, we make some pivotal observations.
\begin{fact} \label{fact:LQGs}
Under Assumption \ref{ass:affine_pg}, the following hold:
\begin{enumerate}[(i)]
\item The PG $ F $ is affine;
\item The PPG map $ h(*,\cdot) =\proj{\ml{Y}(*)}[~\cdot - ~\gamma F(*, \cdot)]$ is PWA;
\item The Jacobian \( \jac h(x, \cdot) \) is well-defined and constant within the partitions \( \{ \tilde{\ml{Y}}_i \}_{i \in \polyset} \), namely, for each $i \in \mathcal N_p$
\begin{equation}
\jac h(x, \cy) = R_i, \quad \forall \cy \text{ such that }  \com(x,\cy) \in \textrm{int } \tilde{\ml{Y}}_i,
\end{equation}
where \( R_i = [R_{i,1},\, R_{i,2}] := [\jac_1 h(x, \cy), \, \jac_2 h(x, \cy)] \).
{\hfill $\square$}
\end{enumerate}
\end{fact}

\begin{figure}
		\centering
		\def\svgwidth{0.95 \linewidth}
       \includegraphics[width=.4\textwidth]{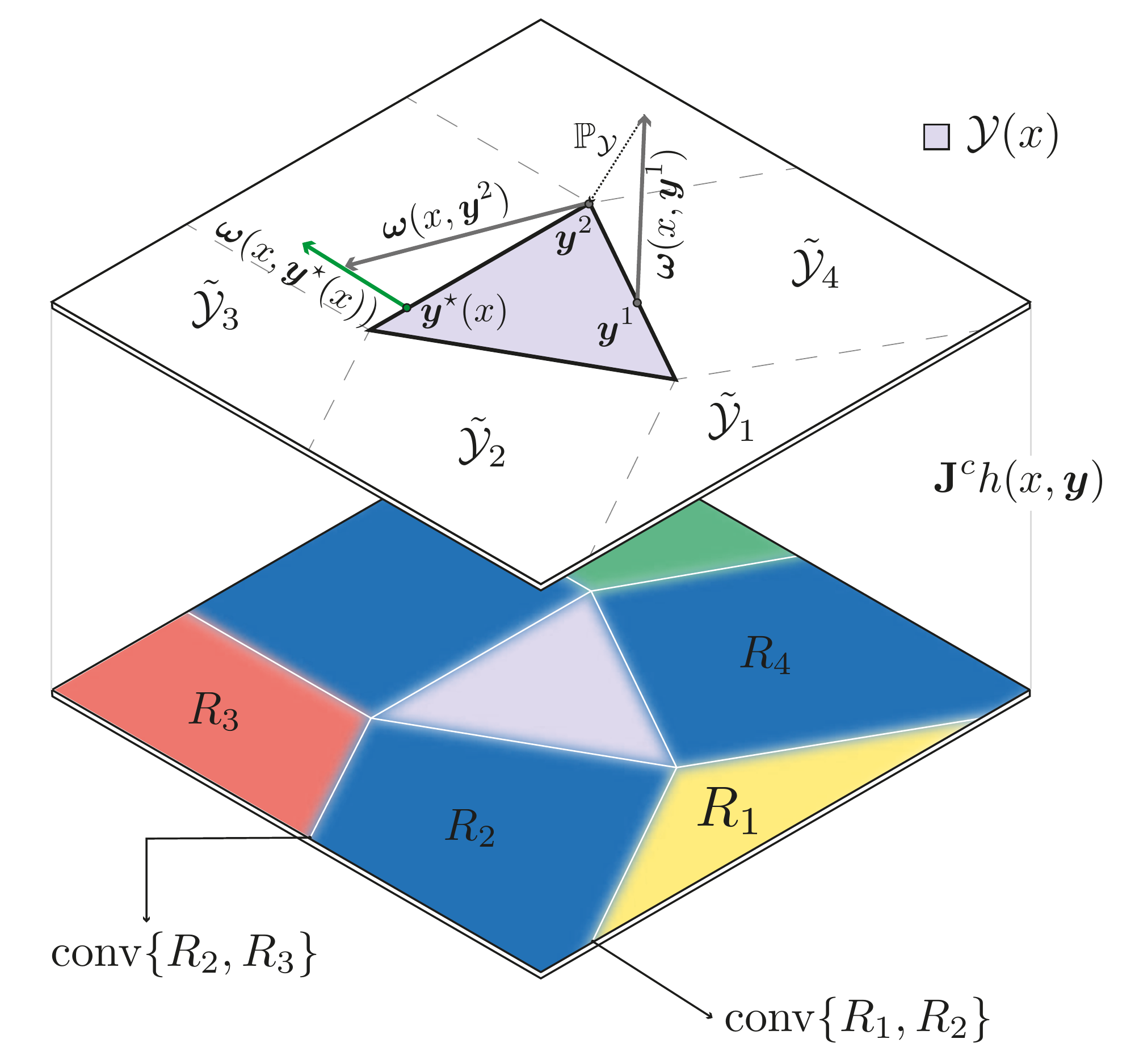}
		\caption{Top: Part of a polyhedral partition
		\( \{\tilde{\ml{Y}}_i \}_{i \in \polyset} \) of \( \setR^2 \).
		Bottom: For LQGs, if \( \com(x,\cy) \) is in the interior of \( \tilde{\ml{Y}}_i\), then \( \clarke h(x, \cy)\) is constant and equal to the standard Jacobian \( \mathbf J h(x, \cy)=R_i \) (coloured areas). 
             If instead \( \com(x,\cy) \) lies on the common boundary between partitions, \( \clarke h (x, \cy) \) 
            is the convex hull of the Jacobians of the active partitions (solid white lines).}
		\label{fig:feasible_sets}
	\end{figure}

With these observations in mind, we can explicitly express the Clarke Jacobian of $ h({x}, \cdot)$ at $\cys({x})$ as
\begin{equation}
\clarke h({x}, \cys({x})) = \conv \{R_i \, | \, i \in \ml{P}({x}) \},
\end{equation}
where \( \ml{P}({x}) \subseteq \polyset \) is the set of active partitions at $\com(x,\cy^\star(x))$, and it is not a singleton if and only if \( \proj{\ml{Y}({x})}[\, \cdot\, ] \) is nondifferentiable at \( \com({ x}, \cys({x})) \). A schematic illustration of this technical discussion is provided in \autoref{fig:feasible_sets}.

The fact that \( \jac h(x, \cdot) \) has a piece-wise constant structure allows for both a theoretical and 
	computational improvement of the sensitivity learning \eqref{eq:s_update_perturbed}.
	To illustrate this, we discuss next the smooth and non-smooth cases separately.

\textbf{Smooth regime:} Assume that \( x \in \ml{X} \) is such that \( \ml{P}(x) \) is a singleton. Whenever \( \tilde{\cy}^\ell \) and \( \cys(x) \) are in the same partition	\(  \tilde{\ml{Y}}_i  \), the nominal sensitivity iteration \eqref{eq:s_update_nominal} and the perturbed one \eqref{eq:s_update_perturbed}
	coincide.
	Hence, convergence is guaranteed by \autoref{prop:s_nominal_convergence}, and the matrices 
        \( \jac_1 h(x, \tilde{\cy}^{\ell+1}), \jac_2 h(x, \tilde{\cy}^{\ell+1}) \) in the sensitivity update remain constant.

\textbf{Nonsmooth regime:} If instead \( x \in \ml{X} \) is such that \( \ml{P}(x) \) is not a singleton, the nominal 
	iteration \eqref{eq:s_update_nominal} is not well-defined, since $h(x,\cdot)$ is not differentiable at $\cy^\star(x)$.
	Nonetheless, we can show that by using any pairs \( [R_{i,1}, R_{i,2}] \in \clarke (x,\cy^\star(x)) \), for some $i \in \ml{P}(x)$, as surrogate partial Jacobians to
	perform \eqref{eq:s_update_perturbed}, 
	the iterates \( \{\tilde{s}^\ell\}_{\ell \in \setN} \) converge to an element of
        the conservative Jacobian\footnote{Interestingly, in general one cannot expect this modified sensitivity learning iteration to converge to an element of the Clarke Jacobian $\clarke \cy^\star(x) $, even when $[R_{i,1}, R_{i,2}]$ is an element of the Clarke Jacobian \( \clarke h(x,\cy^\star(x)) \) \cite{nonsmooth_implicit}.} \( \conserv \cys(x) \) (see section~\ref{subsec:NP} for a formal definition).
	Intuitively, this approach generates a directional derivative of \( \cys(x) \) where the choice of \( R_i \) determines the direction.
             
	\smallskip
	Based on these observations, we propose the inner loop tailored for LQGs in Alg.~\ref{alg:innner_general_setup} (LQG). To improve performance, we replace the Jacobian Update in \autoref{alg:innner_general_setup} by the conditional update in \autoref{alg:nonsmooth_sensitivity_learning_lqgs} that extends the sensitivity learning scheme to nonsmooth regimes.
	\begin{figure}[t]
			\begin{algorithm}[Jacobian Update for LQGs]				\label{alg:nonsmooth_sensitivity_learning_lqgs}
$
\begin{array}{l}
\text{Jacobian conditional update:}\\
\left|
\begin{array}{l}
\text{If } \|\tilde \cy^{\ell+1}- \tilde \cy^\ell \| \geq \sigma: \\[.2em]
\left|
\begin{array}{l l}
S_{1,i} \in \clarke_1 h_i(x, \tilde{\conc{y}}^{\ell+1}), & \text{ 
via \eqref{eq:ppg_jac_1_i}}\\
S_{2,i} \in \clarke_2 h_i(x, \tilde{\conc{y}}^{\ell+1}), & \text{  via \eqref{eq:ppg_jac_2_i}} 
\end{array}
\right.
\\[.3em]
\end{array}
\right.\\[1.5em]
%
\end{array}						~~~~~~~~~~~~~~~~~~~~~~~~~~~~~~~~~~~~~~~~~~~~~~~~~~~~~~~~~~~~~~~~~~~~~~~~~~~~~~~~~
$
				\end{algorithm}
		\end{figure}
		%
            In there, for sufficiently small tolerances, the Jacobian update stops whenever \( \tilde{\cy}^{\ell} \) lies in \( \tilde{\ml{Y}_j} \) for some \( j \in \ml{P}(x) \).
            Then, \( S_1 \) and \( S_2 \) are no longer updated and \( [S_1,\, S_2] = R_j \), as desired.
		Convergence of the inner loop with nonsmooth sensitivity updates is established next.
		\begin{proposition} \label{prop:s_lqg_convergence}
			Under \autoref{ass:affine_pg}, for any \( x \in \ml{X} \) and \( \sigma \) small enough, \( \{\tilde{s}^\ell \}_{\ell \in \setN} \) generated by \autoref{alg:innner_general_setup} with sensitivity updates in \autoref{alg:nonsmooth_sensitivity_learning_lqgs} 
			converges to an element of \( \conserv \cys(x) \).
			\hfill \( \square \)
		\end{proposition}
		
		Now, we exploit the structure of LQGs to derive tighter error bounds for sensitivity learning.
		\begin{lemma} \label{lemma:s_bound_lqg}
			Fix \( x \in \ml{X} \) and let \( \ell' \) be such that \( \com(x, \tilde{\cy}^{\ell'}) \in \tilde{\ml{Y}}_j \) for
			some \( j \in \ml{P}(x) \).
			Under Assumption \ref{ass:affine_pg}, for any \( \ell \geq \ell' \), the iterates \( \tilde{s}^\ell \) of  \autoref{alg:innner_general_setup} with Jacobian updates in 
                \autoref{alg:nonsmooth_sensitivity_learning_lqgs} satisfy 
			\begin{equation} \label{eq:s_bound_lqg}
				\dist \big(\tilde{s}^\ell; \, \conserv \cys(x) \big) 
                    \leq \frac{1}{1 - \eta} \norm{\tilde{s}^{\ell+1} - \tilde{s}^\ell},
			\end{equation}
whenever the tolerance \( \sigma \) is set sufficiently small.
{\hfill $\square$}				
		\end{lemma} 
	
		\subsubsection{Linear-Quadratic Subspace Games (LQSGs)}
		\label{subsubsec:LQGs}
		Finally, we focus on the subclass of LQSGs and exploit their additional structure to further strengthen our sensitivity
		learning results. We begin with some key observations, detailed next.
\begin{fact}
Under Assumptions \ref{ass:affine_pg} and \ref{ass:affine_subspace_con}, the following hold:%
\begin{enumerate}[(i)]%
\item The projection \( \proj{\ml{Y}(x)}[\cdot] \) is affine;
\item The PPG mapping \( h(x,\cdot) \) is affine;
\item The Jacobian mapping \( \jac h(x, \cdot) \) is everywhere well-defined and constant-valued.
{\hfill $\square$}
\end{enumerate}
\end{fact}

	In light of these observations, we conclude that the nominal 
		\eqref{eq:s_update_nominal} and perturbed \eqref{eq:s_update_perturbed} sensitivity iterations coincide.
		Moreover, thanks to the linear dependence of \( \ml{Y}(\cdot) \) on \( x \),
		it can be shown that \( \jac h(\cdot, \cdot) \) is constant with respect to
		both arguments, see Appendix \ref{app:diff_proj}.
		Hence, \( \jac h(\cdot, \cdot) \) need only be computed once, regardless
		of the updates of \( x \) and \( \cy \).
		Finally, recalling that the solution mapping is affine, namely, \( \cys(x) = W x + w \), for some \( W \in \setR^{n \times m}\) and \(w\in \setR^{n} \), we deduce that
		the sensitivity \( \jac \cys(\cdot) \) is constant and equal to \( W \).
		These considerations are formalized in the following proposition.
		\begin{proposition} \label{prop:s_lqgs_convergence}
			Under Assumptions \ref{ass:affine_pg} and \ref{ass:affine_subspace_con},  for any \( x \in \ml{X} \), \( \jac \cys(x) = W \), for some \( W \in \setR^{n \times m}\),  and
			the sequence \( \{\tilde{s}^\ell\}_{\ell \in \setN} \) generated by \eqref{eq:s_update_perturbed}
			converges to \( W \) linearly.	
{			\hfill \( \square \)}
		\end{proposition}
		
	In light of this result, the inner loop for LQSGs reduces to a single
        iteration of \eqref{eq:y_update_distributed} and \eqref{eq:s_update_distributed},
        as shown in \autoref{alg:innner_general_setup} (LQSG).
    This improvement enabled by the stronger error bounds derived in the following lemma.
	
	\begin{lemma} \label{lemma:y_s_lqgs_bound}
		Under Assumptions \ref{ass:affine_pg} and \ref{ass:affine_subspace_con}, there exist constants \( B_{Y1}, B_{Y2}, B_{S} > 0 \) such that the iterates \( \cy^k \) 
		and \( s^k \) of \autoref{alg:general_hypergradient}, with the inner loop in \autoref{alg:innner_general_setup} (LQSG), satisfy
		\begin{align}
		\textstyle
			\text{(i)} & \textstyle
			~ \norm{\cy^{k} - \cys(x^k)} \leq B_{Y1} \eta^{k} 
			+ B_{Y2} \sum_{\ell=0}^{k-1} \eta^{k-\ell} \beta^\ell
			\label{eq:y_vanish_bound} \\
			\text{(ii)} & ~ \norm{s^{k} - W}  \leq B_{S} \eta^{k}   \label{eq:sens_vanish_bound}
		\end{align}
	\end{lemma}
%

	\subsection{Outer Loop Analysis}
	\label{subsec:outer_loop_an}
	\subsubsection{General Games}
    In this section, we prove convergence of \autoref{alg:general_hypergradient} to \textit{composite critical points} of \eqref{eq:upper_level_substi}, 
    namely, some \( x \in \ml{X} \) that satisfies the inclusion
	\begin{equation} \label{eq:composite_critical_points}
		0 \in \conserv \varphi_e(x) + \ncone_{\ml{X}}(x),
	\end{equation}
	where $\ncone_{\ml{X}}$ is the normal cone operator \cite[Def.~6.38]{bauschke2017} of the set $\ml{X}$.
   \blue{Note that, by virtue of \cite[Prop.\ 1]{conservative_definition}, \eqref{eq:composite_critical_points} is a necessary condition
   for \( x \) to be a local minimum of \eqref{eq:upper_level_substi}, which is also referred to as 
   a local Stackelberg equilibrium \cite[Def.\ 1]{local_stack_seeking}.}
    In general, there may be points that are composite critical but are not local minima of \eqref{eq:upper_level_substi}, such as local maxima. While we do not have guarantees that BIG Hype avoids such points, we conjecture that the approximate hypergradient computation helps in that respect, as supported by our numerical experience.

    Our analysis consists of two main steps. In the first step, we rewrite \autoref{alg:general_hypergradient},\blue{ with \( \beta^k = 1\), for all \( k \in \setN \)}, as the iteration
    \begin{equation} \label{eq:leaders_generic_iteration}
        x^{k+1} = x^k + \alpha^k (\xi^k + e^k),
    \end{equation}
    where \( \xi^k \in - \conserv \varphi_e(x^k) - \ncone_{\ml{X}}(x^k) \),
    and \( e^k \) is the hypergradient approximation error. 
    Specifically, the elements of \eqref{eq:leaders_generic_iteration} read as
    \begin{equation}
        \xi^k := - \frac{1}{\alpha^k}(x^k - \proj{\ml{X}}[x^k - \alpha^k \zeta^k] )
    \end{equation}
    where \( \zeta^k \in \conserv \varphi_e(x^k)  \), while the error is given by 
    \begin{align}
	e^k & := \frac{1}{\alpha^k} \big( \proj{\ml{X}} \big[x^k - \alpha^k \widehat{\nabla} \varphi_e^k \big] - \proj{\ml{X}} \big[x^k 
        - \alpha^k \hat{\zeta}^k \big] \big) \label{eq:error_def_general} \\
		\hat{\zeta}^k & \in \argmin_{\zeta \in \conserv \varphi_e(x^k)} \lVert \proj{\ml{X}}[x^k - \alpha^k \widehat{\nabla} 
            \varphi_e^k] - \proj{\ml{X}}[x^k - \alpha^k \zeta]  \rVert.
	\end{align}
	If \( \varphi_e \) is differentiable at \( x^k \), then
        \( \hat{\zeta}^k \!=\!  \{ \nabla \varphi_e(x^k) \} \!=\! \conserv \varphi_e(x^k) \).\\
In the second step, we invoke \cite[Th.\ 3.2]{stochastic_sub_tame} that establishes convergence of \eqref{eq:leaders_generic_iteration} under the following technical conditions:
 	\begin{enumerate}[(i)]
		\item the function \( \varphi_e \) is definable;
		\item the weighted error sequence \( \{\alpha^k e^k\}_{k \in \setN} \) is summable.
	\end{enumerate}
        The first condition follows by definability of \( F \) and \( \varphi \).
        	\begin{lemma} \label{lemma:leader_is_definable}
		The function \( \varphi_e \) is definable.
		\hfill \( \square \)
	\end{lemma}
       The second condition can be enforced under appropriate design of the step sizes and tolerances.
       Before proving it, we introduce an additional technical assumption on the tolerances.
       \begin{assumption}
       \label{ass:tolerance_distance}
       Let
	\( \delta(x) 
        \)
        be the radius of the largest ball whose elements are all included in one of the active partitions, i.e., \( \delta(x) := \max\{ r \in \setRp \, | \, \ball(\com(x, \cys(x); r) \subseteq \cup_{i \in \ml{P}(x)} \tilde{\mathcal{Y}}_i \}
        \). 
   	There exists \( k' \in \setN \) such that \( \sigma^k < \delta(x^k) \), for all \( k > k' \).
   	\( \hfill \square \)
        \end{assumption}
        
         This assumption ensures that, eventually, the sensitivity learning scheme acts on the correct
	piece of the polyhedral partition. In practice, it always holds if, for example, \autoref{alg:general_hypergradient} converges to a point on which \( \varphi_e \) is differentiable, or if the tolerance sequence is chosen to vanish fast enough.

	\begin{lemma} \label{lemma:summable_errors_general}
        Let \( \{ \alpha^k \}_{k \in \setN} \) be nonnegative, nonsummable and square-summable,
      \( \beta^k = 1 \), for all \( k \in \setN \),
        and let \( \{ \sigma^k \}_{k \in \setN} \) be such that \( \sum_{k=0}^{\infty} \alpha^k \sigma^k < \infty \).
		Further, let \autoref{ass:tolerance_distance} hold and assume that \( \ml{P}(x^k) \) is a singleton for all, except  
            finitely many, elements of
		the sequence \( \{x^k\}_{k \in \setN} \) generated by \autoref{alg:general_hypergradient},
		with the inner loop in \autoref{alg:innner_general_setup} (General).
		Then, \( \{\alpha^k e^k\}_{k \in \setN} \),
		with \( e^k \) as in \eqref{eq:error_def_general}, is summable.	
		\hfill \( \square \)
	\end{lemma}
 
	Now, we are ready to prove convergence of \autoref{alg:general_hypergradient}.
	\begin{theorem} \label{th:general_convergence}
        {\color{black} Let Assumption \ref{ass:tolerance_distance} hold.} Assume that \( \ml{P}(x^k) \) is a singleton for all, except  
        finitely many, elements of the sequence \( \{x^k\}_{k \in \setN} \) generated by \autoref{alg:general_hypergradient},
        with \( \{\alpha^k, \beta^k, \sigma^k \}_{k \in \setN} \) as in Lemma \ref{lemma:summable_errors_general} and
        inner loop as in \autoref{alg:innner_general_setup} (General).
        Then, every limit point of \( \{x^k\}_{k \in \setN} \) is composite critical for
        \eqref{eq:upper_level_substi}, and the sequence of function values \( \{ \varphi_e(x^k) \}_{k \in \setN} \)
        converges. \hfill \( \square \)
	\end{theorem}
        
        Assuming the trajectory almost everywhere differentiable is mild in general but essential for conducting meaningful inexact calculations of the hypergradient.
        Next, we present an illustrative example to show that it is not restrictive in practice.%
            	\begin{figure}[b]
    		\centering
    		\def\svgwidth{0.60 \linewidth}
\includegraphics[scale=.6]{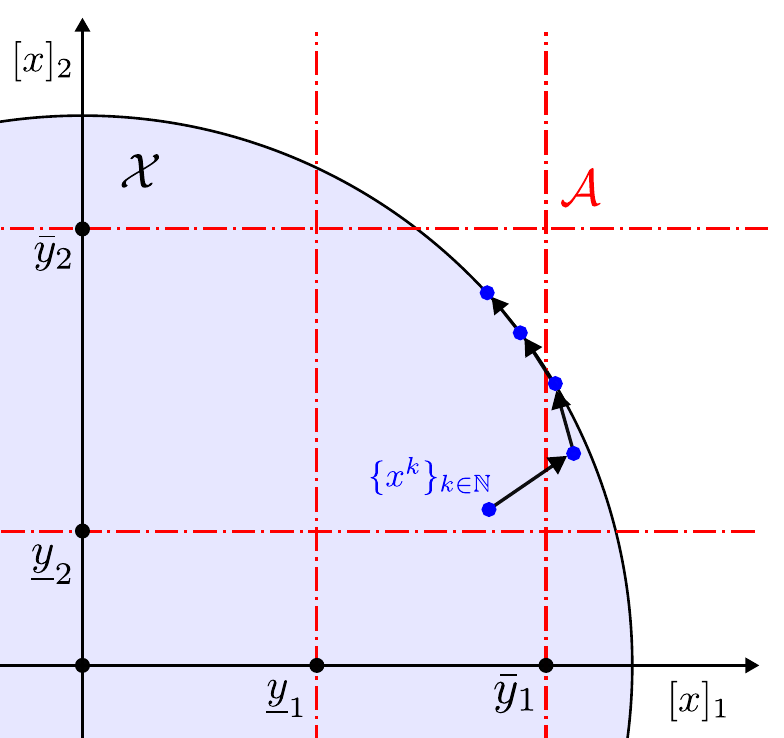}    		
    		\caption{The sets \( \ml{A} \) and \( \ml{X} \) in \autoref{example:differentiable}.}
    		\label{fig:example}
    	\end{figure}
        \begin{example} \label{example:differentiable}
    	Let \( N \!=\! 2 \),
    	\( f_i(x, y_i) := (y_i - [x]_i)^2 \), and \( \mathcal Y_i = [\underline{y}_i, \overline{y}_i] \), for all \( i \in \{ 1,2\} \);
    	moreover let \( \mathcal X = \{\xi \in \mathbb R^2~|~\| \xi \|\leq r \} \) for some \( r > 0 \),
    	and \( \varphi(x, \cy) := - (y_1 + y_2) \).
    	The lower-level solution mapping reads in closed form as \( y_i^{\star}(x) = \max\{\underline{y}_i, \min\{\bar{y}_i, [x]_i\}\} \). Therefore, \( \ml{P}(x) \) is not a singleton if and only if
    	\( x \in \ml{A}:= \{\xi \in \setR^2 \, | \, \exists i \in \{1,2\} \text{ s.t. } [\xi]_i = \bar{y}_i \text{ or } [\xi]_i =                \underline{y}_i\} \).
    	We illustrate \( \ml{A}\) and \( \ml{X} \), and a possible trajectory of \( \{ x^k \}_{k \in \setN} \) in Fig. \ref{fig:example}.
    	The differentiability assumption in Theorem~\ref{th:general_convergence} requires that \( \{x^k\}_{k \in \setN} \) intersects the measure-zero set \( \ml{A} \) on, at most, finitely many points.
		
       \end{example}

	\subsubsection{Linear-Quadratic Games}
        For LQGs, we prove convergence of \autoref{alg:general_hypergradient}
        without assuming almost everywhere differentiable trajectories by leveraging the results in \S~\ref{subsubsec:LQGs}.
	\begin{theorem} \label{th:convergence_lqg}
		Let Assumptions \ref{ass:affine_pg} and \ref{ass:tolerance_distance} hold.
		Then, every limit point of \( \{x^k\}_{k \in \setN} \) generated by \autoref{alg:general_hypergradient},
		with $ \{\alpha^k, \beta^k, \sigma^k \}_{k \in \setN} $ as in Lemma \ref{lemma:summable_errors_general}, the inner loop in \autoref{alg:innner_general_setup} (LQG), and the sensitivity learning step in \autoref{alg:nonsmooth_sensitivity_learning_lqgs}, is composite critical for
		\eqref{eq:upper_level_substi} and the sequence \( \{ \varphi_e(x^k) \}_{k \in \setN} \)
		converges. \hfill \( \square \)
	\end{theorem}

	\subsubsection{Linear-Quadratic Subspace Games}
	Finally, for LQSGs we mitigate the hypergradient approximation error using a Krasnosel'skii--Mann (KM) relaxation step
	with size \( \{ \beta^k \}_{k \in \setN} \).
	Technically, we interpret \autoref{alg:general_hypergradient} (LQSG),
	as an inexact KM iteration \cite[Eqn.\ 5.42]{bauschke2017}
	and invoke \cite[Prop.\ 5.34]{bauschke2017} to establish convergence.
	Similar to the previous analyses, the main technical challenge is showing that the relaxed error
	sequence \( \{ \beta^k \norm{e^k} \}_{k \in \setN} \) is summable, where
	\begin{equation} \label{eq:error_def_lqgs}
		e^k := \proj{\ml{X}}[x^k - \alpha^k \widehat{\nabla} \varphi_e^k] 
		- \proj{\ml{X}}[x^k - \alpha^k \nabla \varphi_e(x^k)].
	\end{equation}

	The desired result is achieved by leveraging the improved error bounds derived in Lemma \ref{lemma:y_s_lqgs_bound}.
	\begin{lemma} \label{lemma:summable_errors_lqgs}
		Under Assumptions \ref{ass:affine_pg} and \ref{ass:affine_subspace_con},
            let $\alpha^k=\alpha  \) for all \( k \in \setN \), for some $\alpha \in (0, L_{J\varphi_e}^{-1}]$ where \( L_{J\varphi_e} \) is 
		the Lipschitz constant\footnote{Lipschitz continuity of \( \nabla \varphi_e \) is directly implied by 
		the fact that \( \cys(\cdot) \) is affine and that \( \nabla_1 \varphi, \nabla_2 \varphi \)
		are Lipschitz.} of \( \nabla \varphi_e \). Further, let \( \{ \beta^k \}_{k \in \setN} \) be nonsummable, square-summable, and
            satisfying \( 0 \leq \beta^{k+1} \leq \beta^k \leq 1 \) for all \( k \in \setN \).
		Then, the error sequence \( \{ \beta^k \norm{e^k} \}_{k \in \setN} \),
		with \( e^k \) as in \eqref{eq:error_def_lqgs} is summable.
		\hfill \( \square \)
	\end{lemma}
	
	Now, we are ready to prove convergence of our scheme.
	\begin{theorem} \label{th:lqgs_convergence}
		Let Assumptions \ref{ass:affine_pg} and \ref{ass:affine_subspace_con} hold,
            and assume that \( \varphi \) is convex.
		Then, the sequence \( \{x^k\}_{k \in \setN} \) generated by \autoref{alg:general_hypergradient},
		  with \( \{ \alpha^k, \beta^k \}_{k \in \setN } \) as in \autoref{lemma:summable_errors_lqgs}, and
            inner loop as in \autoref{alg:innner_general_setup} (LQSG), converges to a solution
		of \eqref{eq:upper_level_substi}.
		\hfill \( \square \)
	\end{theorem}

	The additional convexity assumption on \( \varphi \) allows to show nonexpansiveness
        of the projected hypergradient operator.

	Finally, we note that constant step sizes are generally preferred over vanishing ones as
	they exhibit faster convergence rates.
	For that reason, we provide an extension of \autoref{th:lqgs_convergence} that allows
	employing a constant relaxation sequence \( \{ \beta^k \}_{k \in \setN} \),
	whenever \( \varphi \) is strongly convex, \blue{
        and ensures that \( \{ x^k \}_{k \in \setN} \) converges to the 
        solution of \eqref{eq:upper_level_substi} while \( \{ (\cy^k, s^k) \}_{k \in \setN} \)
        converges to the corresponding NE and sensitivity.}
	\begin{theorem} \label{th:lqgs_convergence_strong}
		Let Assumptions \ref{ass:affine_pg} and \ref{ass:affine_subspace_con} hold, and
            assume that \( \varphi \) is \( \sigma_{\varphi} \)-strongly convex.
		Then, there exists constant \( \bar{\beta} > 0 \) such that if \( \{\beta^k\}_{k \in \setN} \) is a constant
		sequence that is equal to \( \beta \in (0, \bar{\beta}) \),
		the sequence \( \{ (x^k, \cy^k, s^k) \}_{k \in \setN} \) generated
		by \autoref{alg:general_hypergradient}, 
            with \( \{ \alpha^k \}_{k \in \setN} \) as in \autoref{lemma:summable_errors_lqgs},
            and inner loop as in \autoref{alg:innner_general_setup} (LQSG),
		converges to \( (x^{\star}, \cys(x^{\star}), W) \)
		at a linear rate, where \( x^{\star} \) is the unique solution of \eqref{eq:upper_level_substi}.
		\hfill \( \square \)
	\end{theorem}

\section{Numerical Simulations}
\label{sec:numerical_simulations}
       \newcommand{\tintervals}{\Lambda}
	In this section, we numerically investigate the convergence and scalability properties
        of \autoref{alg:general_hypergradient} by deploying it on a
	simplified version of the demand-response model in \cite{drm_miqp}.

	\subsection{Demand-Response Model}
	We consider a Distribution System Operator (DSO), serving as the leader, selling electricity to a group \( \agents = \{1, \ldots, N\} \) of smart buildings, playing the role of followers.
	The DSO faces a day-ahead scheduling problem, where the day is discretized
	in \( \tintervals \) intervals of length \( \Delta \tau \)
	indexed by \( \ml{T} := \{1, \ldots, \tintervals\} \),
	and seeks to maximize its revenue by regulating the energy price.
	In turn, the smart buildings aim at covering their daily electricity needs at the minimum cost.
	\subsubsection{DSO}
	We denote by \( p_i \in \setRp^\tintervals \) the energy purchased by building \( i \in \agents \)
	throughout \( \ml{T} \), and by \( \bar{p} : = \sum_{i \in \agents} p_i \) the aggregate purchased energy.
	The DSO designs a pricing map \( h : \setR^\tintervals \to \setRp^\tintervals \) that
	determines the cost of electricity at each time period \( \tau \in \ml{T} \).
	We consider the affine price map 
	\( h ( \bar{p} ) := C_1 \bar{p}+ c_0 \),
	where \( c_0 \in [\underline{c}_0, \overline{c}_0]^\tintervals \) is the baseline price, 
	and  \( C_1 := \diag(c_1) \) with \( c_1 \in [\underline{c}_1, \overline{c}_1]^\tintervals \) representing the marginal price.
	The average value of \( c_0 \) is upper-bounded by \( (\underline{c}_0 + \overline{c}_0)/2 \),
	and a similar constraint is imposed on \( c_1 \).
	The DSO modulates prices via \( c_0 \) and \( c_1 \) with the goal of minimizing
	\begin{equation}
			\varphi(c_0, c_1, \bar{p}) = - h(\bar{p})^{\top} \bar{p},
	\end{equation}
where ${h(\bar{p})}^{\top} \bar{p}$ is the total revenue of the DSO over $\mathcal T$.


	\subsubsection{Smart Buildings}
	Each building incurs a cost of the form
	\begin{equation} \label{eq:prosumer_cost}
		f_i := h(\bar{p})^{\top} p_i + \lambda_b \big(\norm{p_i^{\text{C}}}^2 + \norm{p_i^{\text{DC}}}^2 \big),
	\end{equation}
        where the first term is the electricity procurement cost, while the second term models the battery degradation cost. In there,  \( p_i^{\text{C}}, p_i^{\text{DC}} \in \setRp^\tintervals  \) are the charging and discharging power of the battery, respectively, and \( \lambda_b > 0\) is a constant. \blue{The objective \eqref{eq:prosumer_cost} is solely influenced by local variables and by the aggregate demand $\bar p$, thus preventing smart buildings from requiring access to the specific consumption profiles of others.}
The local variables are subject to a series of affine constraints that model the energy demand, the dynamics and charging/discharging limits of the storage device \cite[Eqn.\ 9]{drm_miqp}.

	\subsubsection{Aggregate Load Constraints} \label{subsec:agg_load}
	The aggregate energy consumption of the community must satisfy
	\( \bar{p} \leq g \),
	where \( g \in \setRp^\tintervals \) is the capacity of the grid.
	Following \autoref{subsec:coupling_constraints}, we decompose this coupling constraint by using agent-specific auxiliary variables, coupled through the leader's constraints.
	To this end, we introduce the variable \( \hat{\theta}\in \setRp^N \) and we reformulate the coupling constraint \( \bar{p} \leq g \) via the local constraints
	\( p_i \leq \hat{\theta}_i g \), with $i=1, \ldots,N$, and the leader's constraint \( \ones_N^{\top} \hat{\theta} = 1 \).
	Intuitively, \( \hat{\theta}_i \) corresponds to the percentage
	of grid capacity that smart building \( i \in \agents \) is allowed to access.
	We stress that \( \hat{\theta}\) is a design variable that the leader optimizes to determine the allocation of the shared 
	resource \( g \). At first sight, introducing \( \hat{\theta} \) seems to increase the dimension
		of the leader decision vector and, in turn, also the computational burden of its updates.
		We observe, however, that \( \hat{\theta} \) is decoupled from
		the remaining decision variables and its feasible set is a unit simplex \( \ones_N^{\top} \hat{\theta} = 1 \), for which projection can be performed at virtually zero computational cost \cite{condat2016simplex}.

	\subsection{Numerical Results}
			\begin{figure*}[t]
		\begin{subfigure}{0.33 \textwidth}
			\centering
			\includegraphics[width=\linewidth]{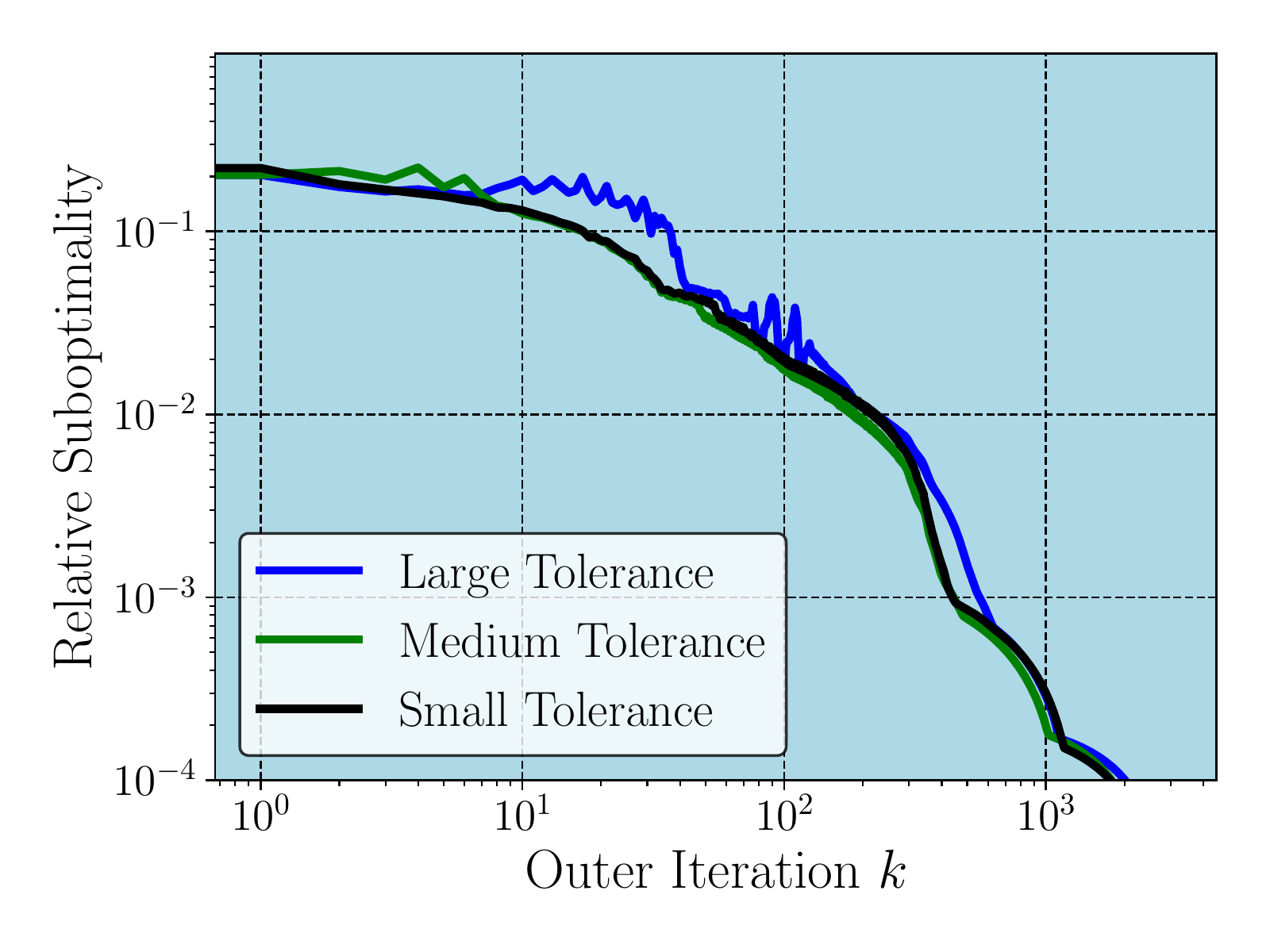}
			\caption{}
			\label{subfig:rel_subopt_outer}
		\end{subfigure}
		\begin{subfigure}{0.33 \textwidth}
			\centering
			\includegraphics[width=\linewidth]{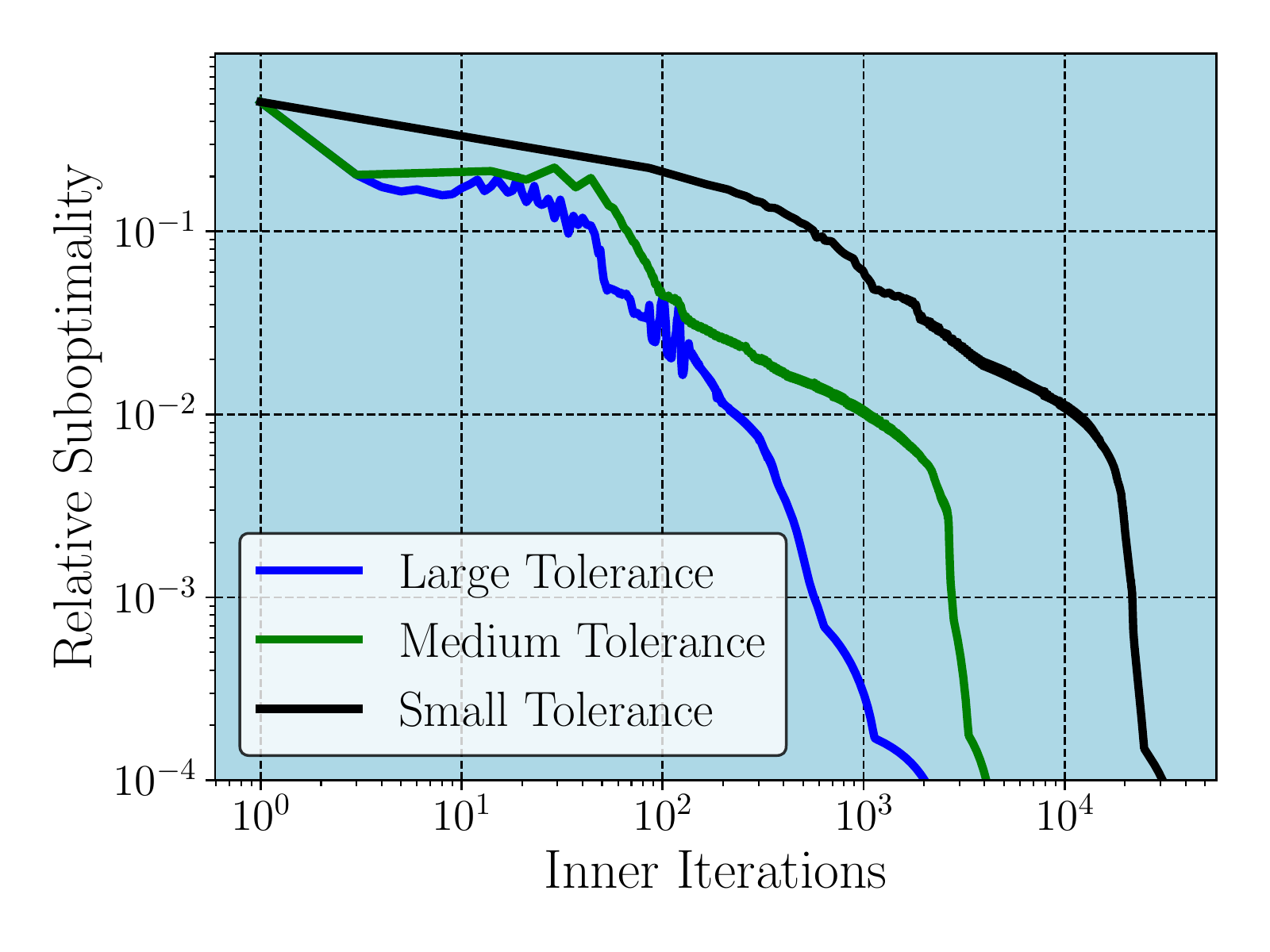}
			\caption{}
			\label{subfig:rel_subopt_inner}
		\end{subfigure}
		\begin{subfigure}{0.33 \textwidth}
			\centering
			\includegraphics[width=\linewidth]{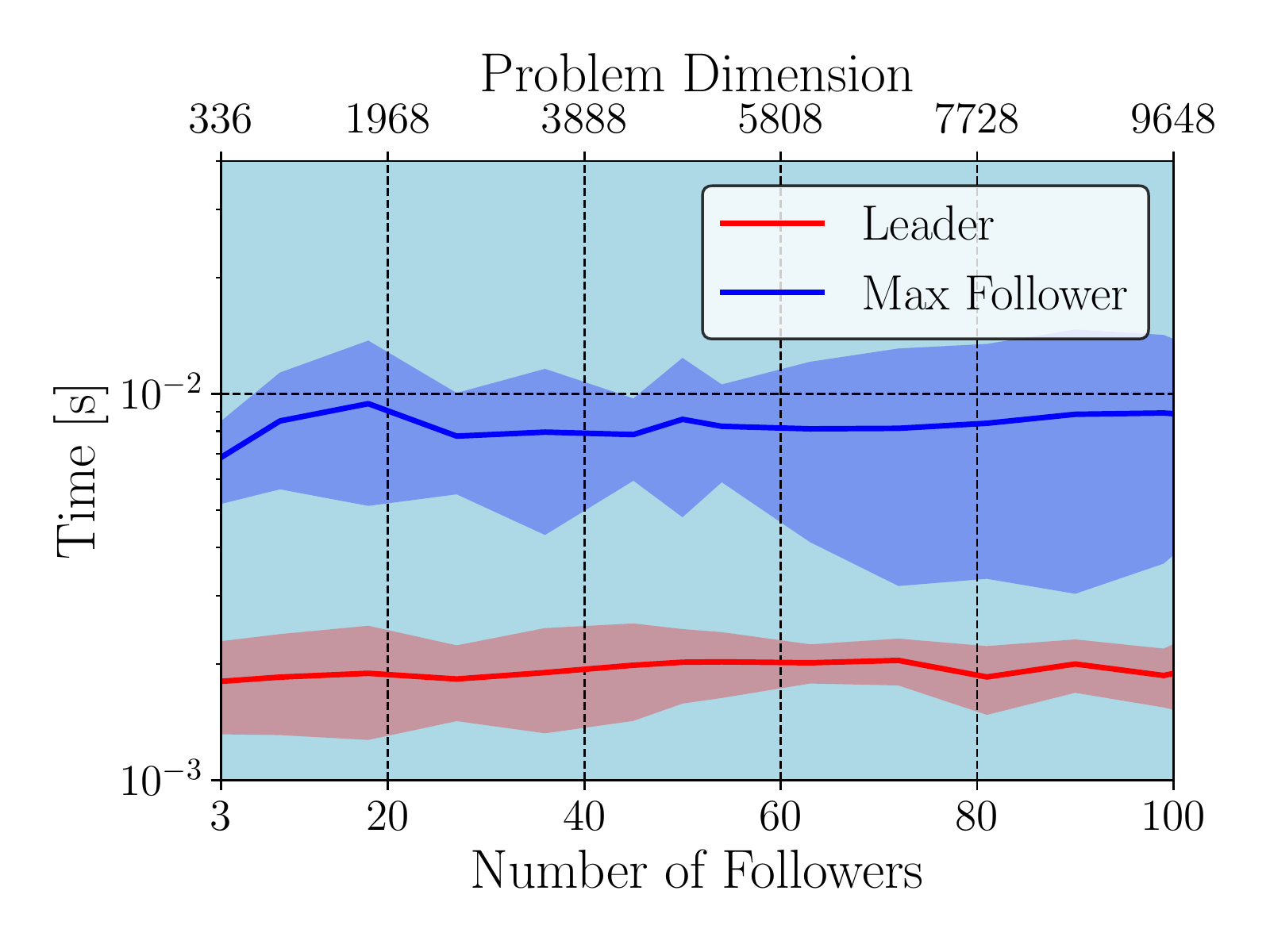}
			\caption{}
			\label{subfig:separate_times}
		\end{subfigure}
		\hfill
		\caption{(a), (b) Relative suboptimality as a function of outer and total inner-loop iterations, respectively. 
		(c) Computational time for inner and outer loop iterations versus number of followers and total problem dimension.}
		\label{fig:convergence_behavior_test}
	\end{figure*}
	\subsubsection{Implementation Details}
	We use OSQP \cite{osqp} to compute polyhedral projections,
	while for their Jacobians we solve a series of linear systems
	as discussed in Appendix \ref{app:diff_proj}.
	All simulations are implemented in Python and run on a MacBook Air with an M1 chip and 16GB RAM.
	\blue{The code of the simulations and a general Python implementation of BIG Hype is available online at \nolinkurl{https://github.com/p-grontas/BIG_hype_algorithm}.}
	\subsubsection{Convergence Behavior}

\blue{			We consider an energy community comprising \( N = 9 \) smart buildings and investigate the convergence behavior of \autoref{alg:general_hypergradient} with respect to different choices of the
			tolerance sequence.
			In practice, we select two separate tolerance
			sequences \( \{\sigma_{\cy}^k\}_{k \in \setN} \)
			and \( \{\sigma_{s}^k\}_{k \in \setN} \) to independently fine-tune the equilibrium and sensitivity errors, respectively, while maintaining theoretical guarantees.
		We use tolerance sequences of the form \( a_{\cy} (k+1)^{-0.51}\), \(a_{s} (k+1)^{-0.51} \), and study three choices of the parameters \( a_{\cy}, a_{s} \). Specifically,
		we consider small tolerances \( (a_{\cy}, a_{s}) = (0.002, 25) \), which approximate running the inner
		loop to full convergence; large tolerances \( (a_{\cy}, a_{s}) = (0.1, 500) \), which typically translates into few inner loop iterations and inaccurate hypergradients; 
		and medium tolerances \( (a_{\cy}, a_{s}) = (0.02, 50) \).
		The step size of the hypergradient is set as 
		\( \alpha^k = 3  \cdot 10^{-6} (k+1)^{-0.51} \).}

	We run \autoref{alg:general_hypergradient}, under all three tolerance choices
	and compare the resulting relative suboptimality, defined as \( (\varphi_e(x^k) - \varphi_e^{\star})/\abs{\varphi_e^{\star}} \)
	where \( \varphi_e^{\star} \) is the lowest objective value obtained over all runs.
	The relative suboptimality as a function of outer loop iterations \( k \) is shown in \autoref{subfig:rel_subopt_outer}.
	We observe that medium and small tolerances result in slightly faster convergence,
	whereas inaccurate hypergradients give rise to oscillations.
\autoref{subfig:rel_subopt_inner} illustrates the relative suboptimality with the respect to the number of inner loop iterations that better reflects the cost of computation and communication. We observe that large tolerances are preferable since they allow for fewer inner iterations without compromising convergence.

        \subsubsection{Comparison with Literature}
        We compare the computational cost and solution quality of \autoref{alg:general_hypergradient}
        to a popular solution approach, namely, mixed-integer programming.
        To cast \eqref{eq:upper_level_expanded} as a mixed-integer program (MIP), 
        we use the big-M reformulation \cite{FortunyAmat1981ARA}, and we use Gurobi to solve the resulting problem,
        which is a MIP with quadratic objective and nonconvex quadratic constraints.
        The solver is terminated once a relative suboptimality of \( 5 \% \) is attained.
        We compare BIG Hype and MIP in terms of computation time and solution quality, for different problem sizes.
        For BIG Hype, we consider the computation time for both \( 1 \) CPU and \( N \) CPUs, 
        corresponding to a serial and a distributed implementation, respectively.
        The computational time of the inner loop for \( N \) CPUs is calculated as the maximum time over all agents.
        The termination condition for \autoref{alg:general_hypergradient} is
        \( \abs{\varphi_e(x^{k+1}) - \varphi_e(x^{k})}/\abs{\varphi_e(x^{k})} \leq 10^{-5} \).
        We quantify the solution quality as
        \blue{\( (\varphi_{\textrm{all}}^{\star} - \varphi_{\textrm{alg}}^{\star}) / \abs{\varphi_{\textrm{all}}^{\star}} \)},
        where \( \varphi_{\textrm{alg}}^{\star} \) is the best objective value attained by each algorithm at 
        convergence or time out, 
        and \( \varphi_{\textrm{all}}^{\star} \) is the best objective over both algorithms.
        All results are summarized in \autoref{table:comparison_literature}.        
        We note that BIG Hype consistently achieves faster computation times, even with small problem sizes. Notably, when $N=1$, both methods yield virtually identical solutions. 
        However, already for \( N = 2 \), the MIP approach can only produce a significantly suboptimal
        solution within the allotted time frame.
        Overall, BIG Hype demonstrates rapid convergence without compromising the quality of the solution.
                	\begin{table}[htbp]
        	\centering
        	\begin{tabular}{r||c|c|c}
        		\( (N, \tintervals) \)		 & MIP                  & BIG Hype (\( 1 \) CPU)  & BIG Hype (\( N \) CPUs)  \\ \hline
        		\(  (1, 8) \)  & 4hr (\( 0\% \))                     & 5s (\( -0.1\% \))             & 5s (\( -0.1\% \))                    \\ 
        		\( (1, 24) \)  & 6hr (\( -0.4\% \))  		& 14s (\( 0\% \))			  & 14s (\( 0\% \))                       \\
        		\( (2, 24) \)  & 6hr (\( -31\% \)) 		& 12s (\( 0\% \))			  & 7s (\( 0\% \))                      
        	\end{tabular}
        	\caption{Computation times and solution quality (in parentheses) for MIP and BIG Hype across various problem sizes.
}
        	\label{table:comparison_literature}
        \end{table} 

	\subsubsection{Computational Cost}
	Now, we investigate the computational cost and the scalability of \autoref{alg:general_hypergradient}
	with respect to the number of smart buildings.
	We consider \( N \) up to 100 followers, and determine
	the computation time required for each inner and outer loop update of \autoref{alg:general_hypergradient}.
	For the inner loop, we approximate the performance of a distributed implementation by considering the
	maximum computational cost among all followers.
	Our results are presented in \autoref{subfig:separate_times}. 
        Notice that the computational cost remains
	constant due to the aggregative nature of the demand response model.
	The computational burden of the followers is considerably larger, as
	due to the projection operation (and its Jacobian) onto a higher-dimensional feasible set with respect to the leader.
	
\section{Conclusion}
\label{sec:conclusion}
Hypergradient descent can be effectively deployed to find local equilibria of large-scale Stackelberg games. Not only this approach is scalable, thus making it suitable for realistic applications but also reliable, as convergence is mathematically guaranteed for appropriate choices of its hyperparameters. Additionally, for liner-quadratic games with polyhedral and affine constraints, respectively, hypergradient descent can be refined to yield both theoretical and numerical improvements.

\blue{
This approach opens a number of future research directions.
On the theoretical side, it would be interesting to investigate the use of constant step sizes that would enable faster convergence, and pave the way for an online implementation.
On the practical side, BIG Hype could be deployed to solve sophisticated and large-dimensional hierarchical problems in different application domains, such as traffic routing and energy trading, which proved to be computationally challenging.}

\appendices

\section*{Appendix}
\renewcommand\thesubsection{\Alph{subsection}}
	\subsection{Differentiating Through a Projection} \label{app:diff_proj}
	To compute the Jacobian of the polyhedral projection 
\begin{align*}
\proj{\ml{Y}_i(x)}[y] = 
\left\{
\begin{array}{r l} 
\argmin_{z } &\frac{1}{2}z^\top I z - y^\top z\\
\textrm{subject to} & A_i z  \leq b_i + G_i x \\
& C_i y = d_i + H_i x
\end{array}
\right.
\end{align*}
we first write the KKT conditions of the associated quadratic program at a primal-dual solution \( (z^{\star}, \lambda^{\star}, \nu^{\star}) \):
	\begin{subequations}
	\label{eq:KKTsys}
		\begin{align}
			z^{\star} - y + A_i^{\top} \lambda^{\star} + C_i^{\top} \nu^{\star} & = 0, \\ 
			0 \leq \lambda^{\star} \perp - (A_i z^{\star} - b_i - G_i x ) & \geq 0,
			\label{eq:proj_complementarity}\\
			C_i z^{\star} - d_i - H_i x & = 0.
		\end{align}
	\end{subequations}
	Then, we differentiate \eqref{eq:KKTsys} with respect to \( (z^{\star}, \lambda^{\star}, \nu^{\star}) \), as explained in \cite{optnet_arxiv}, yielding the linear system
	\begin{equation} \label{eq:proj_kkt_diff}
		\begin{gathered}
			\begin{bmatrix}
				I & A_i^{\top} & C_i^{\top} \\
				\diag(\lambda^{\star}) A_i & \diag(b_i \!+\! G_ix \!-\! A_i z^{\star}) & 0 \\
				C_i & 0 & 0
			\end{bmatrix}
			\begin{bmatrix}
				dz \\ d \lambda \\ d \nu
			\end{bmatrix} = \\
			\begin{bmatrix}
				dy \\
				\diag(\lambda^{\star}) G_i dx \\
				H_i dx
			\end{bmatrix}.
		\end{gathered}
	\end{equation}
Define $g_i(x,y) := \proj{\ml{Y}_i(x)}[y]$ for ease of notation. To compute the Jacobian \( \jac_1 g_i(x, y) \) it suffices to substitute \( dx = I \) and \( d y  = 0_{n_i\times n_i}\), solve the system \eqref{eq:proj_kkt_diff}, and extract the $dz$-component of the solution.	Similarly, to obtain \( \jac_2 g_i(x, \cy) \) we substitute \( dx = 0_{m \times m} \) and \( d y  = I\), solve \eqref{eq:proj_kkt_diff}, and take first solution component. Notably, \eqref{eq:proj_kkt_diff} always admits a solution \cite[App.\ C.1]{optnet_arxiv}.	Further, if strict complementarity holds for \eqref{eq:proj_complementarity}, namely, exactly one of \( [\lambda^{\star}]_{j} \) and \( [G_i x + b_i - A_i z^{\star}]_{j} \)
	is zero for all components \( j = 1, \ldots, p_i \),
	then the solution is unique.
	Moreover, if we assume strict complementarity and consider only the set \( \ml{E} \) of active inequality
	constraints and non-zero \( [\lambda^{\star}]_j \),	then we can show that the Jacobian of the projection is independent of \( x, y \)
	for fixed \( \ml{E} \).
	
	\subsubsection{Proof of \autoref{lemma:s_perturbed_bound}} 
	\label{proof:s_perturbed_bound}
	To establish the claim, we first show that \( \tilde{\cy}^\ell \) lies on the same
	partition of \( \{\tilde{\ml{Y}}_i \}_{i \in \polyset} \) as \( \cys(x) \), within
	a finite number of inner iterations.
	Then, based on this, we prove the desired bound.	
        By virtue of \autoref{fact:diff_equivalence}, it holds true that
	\( \com(x, \cys(x)) \in \inter \tilde{\mathcal{Y}}_i \)
	for some \( i \in \polyset \).
	Therefore, there exists \( \epsilon > 0 \) such that
	\( \ball(\com(x, \cys(x));\, \epsilon) \subseteq \inter \tilde{\mathcal{Y}}_i \).
	Hence, by continuity of \( \com(x, \cdot) \) there exists \( \delta > 0  \) such 
	that \( \cy \in \ball(\cys(x);\, \delta) \) implies
	\( \com(x, \cy) \in \ball(\com(x, \cys(x));\, \epsilon) \), and since
	\( \tilde{\cy}^\ell \to \cys(x) \) there exists \( \bar{\ell} \in \setN \)
	that guarantees \( \tilde{\cy}^\ell \in \ball(\cys(x), \delta) \) for all 
	\( \ell \geq \bar{\ell} \).
	Then, \( h(x, \cdot) \) is continuously differentiable
	and, crucially, its partial Jacobian
	at \( \tilde{\cy}^\ell \) and \( \cys(x) \) can be evaluated by the same expression
	since the projection is given by the same affine mapping.
	
	Now, for \( \ell \geq \bar{\ell} \) the following bound holds:
	\begin{align*}
		& \norm{\jac_1 h(x, \tilde{\cy}^\ell) - \jac_1 h(x, \cys(x))} \\
		& \leq \norm{\jac_1 g(x,\com(x, \tilde{\cy}^\ell)) - \jac_1 g(x, \com(x, \cys(x)))} \\
		& + \gamma \big\lVert\jac_2 g(x, \com(x, \tilde{\cy}^\ell)) \jac_1 F(x, \tilde{\cy}^\ell) \big. \\
		& ~~~ - \big. \jac_2 g(x, \com(x, \cys(x)) \jac_1 F(x, \cys(x)) \big\rVert
	\end{align*}
	where we used the definition of \( \com \), and substituted \eqref{eq:ppg_jac_1_i}.
	For the first term, we note that \( \com(x, \tilde{\cy}^\ell) \) and \( \com(x, \cys(x)) \)
	lie on the same partition of \( \{ \tilde{\ml{Y}}_i \}_{i \in \polyset} \) and, thus,
	\( \jac_1 g(x,\com(x, \tilde{\cy}^\ell)) = \jac_1 g(x, \com(x, \cys(x))) \) as shown
	in Appendix \ref{app:diff_proj}; hence, the first term vanishes.
	For the second term, we have that \( \jac_2 g(x, \cdot) \) is constant on each
	partition of \( \{ \tilde{\ml{Y}}_i \}_{i \in \polyset} \), thus we have
	\begin{align*}
		& \gamma \big\lVert\jac_2 g(x, \com(x, \cys(x)) \big(\jac_1 F(x,\tilde{\cy}^\ell) - \jac_1 F(x,\cys(x)) \big) \big\rVert \\
		& \overset{(a)}{\leq} \gamma \norm{\jac_2 g(x, \com(x, \cys(x))} \norm{\jac_1 F(x,\tilde{\cy}^\ell) - \jac_1 F(x,\cys(x))} \\
		& \overset{(b)}{\leq} \gamma L_{JF1} \norm{\tilde{\cy}^\ell - \cys(x)},
	\end{align*}
	where (a) holds by submultiplicativity of matrix norms and (b) by \autoref{ass:followers_pg_general}
	and the fact the \( \norm{\jac_2 g(x, \cdot)} \leq 1 \) since the projection is nonexpansive.
	
	For the second bound, we have
	\begin{align*}
		& \norm{\jac_2 h(x, \tilde{\cy}^\ell) - \jac_2 h(x, \cys(x))} \\
		& \overset{(a)}{=} \big\lVert \jac_2 g(x,\com(x, \tilde{\cy}^\ell)) (I - \gamma \jac_2 F(x, \tilde{\cy}^{\ell})) \big. \\
		&	~~ \big. - \jac_2 g(x, \com(x, \cys(x))) (I - \gamma \jac_2 F(x, \conc{y}^{\star})) \big\rVert \\
		& \overset{(b)}{\leq} \gamma \norm{\jac_2 g(x, \com(x, \cys(x))} \norm{\jac_2 F(x, \tilde{\cy}^{\ell}) - \jac_2 F(x, \conc{y}^{\star})} \\
		& \overset{(c)}{\leq} \gamma L_{JF2} \norm{\tilde{\cy}^{\ell} - \conc{y}^{\star}}
	\end{align*}
	where (a) holds by \eqref{eq:ppg_jac_2_i}, (b) since \(  \jac_2 g(x, \cdot) \) is constant on each
	\( \tilde{\ml{Y}}_i \) and (c) by \autoref{ass:followers_pg_general}.
	\hfill \( \blacksquare \)

\subsection{Proofs of Section~\ref{subsec:inner_loop_an}}
        \subsubsection{Proof of \autoref{prop:s_nominal_convergence}}
        \label{proof:s_nominal_convergence}
			By \autoref{fact:diff_equivalence}, we have that \( h(x, \cdot) = \proj{\ml{Y}(x)}[~\cdot - ~\gamma F(x, \cdot)] \) is continuously differentiable at \( \cys(x) \).
                Moreover, \( h(x, \cdot) \) is \( \eta \)-Lipschitz, for any \( x \in \ml{X} \), and
		      hence, by \cite[Th 2.17(i)]{nonsmooth_approach},
			\( \norm{\jac_2 h(x, \cys(x))}_F \leq \eta < 1 \), where \( \norm{\cdot}_F \) 
			is the Frobenius norm.
			Thus, we have \( \norm{\jac_2 h(x, \cys(x))} < 1 \) since \( \norm{\cdot} < \norm{\cdot}_F \).
			Next, notice that \( \cys(x) \) is the unique solution of the equation
			\( \conc{y} - h(x, \cy) = 0 \)
			where \( h(x, \cy) \) is continuously differentiable and
			\( I - \jac_2 h(x,\cy) \) is non-singular.
			Invoking the implicit function theorem \cite[Th. 1B.1]{implicit_function}, we deduce that 
			\( \cys(\cdot) \) is continuously differentiable at \( x \), and the sensitivity derivative satisfies
			\begin{equation} \label{eq:sensitivity_explicit_formula}
					\jac \cys(x) = (I - \jac_2 h(x,\cys(x)))^{-1} \jac_1 h(x,\cys(x)).
			\end{equation}
		
			Note that \( \norm{\jac_2 h(x, \cys(x))} \!<\! 1 \) implies that \eqref{eq:s_update_nominal} is contractive
			and converges to its unique fixed point.
			We conclude by observing that the fixed point of \eqref{eq:s_update_nominal} satisfies
			\eqref{eq:sensitivity_explicit_formula}.
   {\hfill $\blacksquare$}

\subsubsection{Proof of \autoref{lemma:s_bound_general}}
	\label{proof:s_bound_general}
	By the triangle inequality 
	\begin{equation*}
		\norm{\tilde{s}^\ell - \jac \cys(x)} \leq \norm{\tilde{s}^\ell - \hat{s}^\ell} + \norm{\hat{s}^\ell - \jac \cys(x)},
	\end{equation*}
	where \( \hat{s} \) is defined in \eqref{eq:s_update_nominal}.
        We will prove the claim by individually bounding the terms \( \norm{\tilde{s}^\ell - \hat{s}^\ell} \) 
	and \( \norm{\hat{s}^\ell - \jac \cys(x)} \).
	
	For the sake of brevity, we let 
	\( L^\ell := \jac_1 h(x, \tilde{\cy}^\ell), L^{\star} := \jac_1 h(x, \cys(x)),
	M^\ell := \jac_2 h(x, \tilde{\cy}^\ell),  M^{\star} := \jac_2 h(x, \cys(x)) \).
	Then, for the first term we have
	\begin{align}
		&\norm{\tilde{s}^{\ell+1} - \hat{s}^{\ell+1}} 
		 \overset{(a)}{=} \norm{M^\ell \tilde{s}^\ell + L^\ell - M^{\star} \hat{s}^\ell - L^{\star}} \notag \\
		& \overset{(b)}{\leq} \norm{M^\ell} \norm{\tilde{s}^\ell - \hat{s}^\ell} + \norm{M^\ell - M^{\star}} \norm{\hat{s}^\ell}
		+ \norm{L^{\ell} - L^{\star}} \notag \\
		& \overset{(c)}{\leq} \eta \norm{\tilde{s}^\ell - \hat{s}^\ell} 
		+ \gamma L_{JF2} \norm{\tilde{\cy}^{\ell} - \cys(x)} \norm{\hat{s}^\ell} \notag \\  
		& ~~~~~ + \gamma L_{JF1} \norm{\conc{y}^{\ell} - \conc{y}^{\star}(x)} \notag \\ 
		& \overset{(d)}{\leq} \frac{\gamma L_{JF2} \norm{\hat{s}^\ell} + \gamma L_{JF1}}{1 - \eta} \norm{\tilde{\cy}^{\ell+1}- \tilde{\cy}^{\ell}}
		+ \eta \norm{\tilde{s}^\ell - \hat{s}^\ell}, \label{eq:s_bound_nom_pert}
	\end{align}
	where (a) follows from the update rules \eqref{eq:s_update_nominal} and \eqref{eq:s_update_perturbed},
	(b) by adding and subtracting \( M^\ell \hat{s}^\ell \),
	(c) follows from \autoref{lemma:s_perturbed_bound} and the fact that \( \norm{\jac_2 h(x, \tilde{\cy}^\ell)} \leq \eta \),
	and (d) by \eqref{eq:y_bound}.
	Recall that \( \{\hat{s}^\ell\}_{\ell \in \setN} \) is generated by a contraction 
        (see the proof of \autoref{prop:s_nominal_convergence} in Appendix \ref{proof:s_nominal_convergence}), 
        hence there exists \( B_{\hat{s}} > 0 \) such that \( \norm{\hat{s}^\ell} \leq B_{\hat{s}} \)
	for all \( \ell \in \setN \).
	We let \( B_{ps} := (\gamma L_{JF2} B_{\hat{s}} + \gamma L_{JF1})/(1 - \eta) \)
	and telescoping \eqref{eq:s_bound_nom_pert} yields:
	\begin{equation*}
		\norm{\tilde{s}^{\ell+1} - \hat{s}^{\ell+1}} \leq 
		B_{ps} \sum_{j=0}^{\ell} \eta^{\ell-j} \norm{\tilde{\cy}^{j+1} - \tilde{\cy}^{j}}.
	\end{equation*}

	For the second term, by contractiveness of \eqref{eq:s_update_nominal}
	\begin{equation*}
		\norm{\hat{s}^\ell - \jac \cys(x)} \leq \eta^\ell \norm{\hat{s}^0 - \jac \cy^{\star}(x)}.
	\end{equation*}
        Next, recall that \( \cys(\cdot) \) is locally Lipschitz \cite[Th.\ 2.1]{dafermos_sensitivity}
	on the compact set \( \ml{X} \) and, therefore, \( L_S \)-Lipschitz for some \( L_S > 0 \).
	By virtue of \cite[Th.\ 2.17(i)]{nonsmooth_approach}, \( \norm{\jac \cys(x)} \leq L_S \),
	for any \( x \in \ml{X} \).
	Thus, there exists a constant \( B_{ns} \) such that \( \norm{\hat{s}^0 - \jac \cy^{\star}(x)} \leq B_{ns} \),
	hence concluding the proof.
	\hfill \( \blacksquare \)
	
	\subsubsection{Proof of \autoref{prop:s_general_convergence}}
	\label{proof:s_general_convergence}
	To show the claim, we recall that \eqref{eq:y_update} is a contraction
	and, thus,
	\( \norm{\tilde{\cy}^{\ell+1} - \tilde{\cy}^\ell} \leq \eta^\ell \norm{\tilde{\cy}^{1} - \tilde{\cy}^0} \).
	We combine this with \eqref{eq:s_bound_general} in \autoref{lemma:s_bound_general} to derive the bound
	\begin{align*}
		\norm{\tilde{s}^\ell - \jac \cys(x)} & \leq 
		B_{ps} \norm{\tilde{\cy}^{1} - \tilde{\cy}^0} \sum_{j=0}^{\ell - 1} \eta^{\ell - 1 - j} \eta^j +
		B_{ns} \eta^\ell \\
		& = B_{ps} \norm{\tilde{\cy}^{1} - \tilde{\cy}^0} \eta^{\ell - 1} \ell +
		B_{ns} \eta^\ell.
	\end{align*}
	By letting \( \ell \to \infty \), the right-hand side vanishes.
	\hfill \( \blacksquare \)

	\subsubsection{Proof of \autoref{prop:s_lqg_convergence}}
	\label{proof:s_lqg_convergence}
	Initially, we will employ \cite[Th.\ 5]{nonsmooth_implicit} to characterize the
	conservative Jacobian \( \conserv \cys \) in terms of \( \clarke h \).
	To do so, we define the mapping \( f(x, \cy) := \cy - h(x, \cy) \) that
	satisfies \( f(x, \cys(x)) = 0 \) and \( \cys(\cdot) \) is locally
	Lipschitz by virtue of \cite[Th.\ 2.1]{dafermos_sensitivity}.
        Note that PWA functions are semialgebraic by definition,
        and, hence, definable \cite{nonsmooth_implicit}.
	Therefore, \( f \) is definable because
	\( \proj{\mathcal{Y}(x)}[\cdot]\) is PWA, and the same holds
	for \( h(x, \cdot) \).
	We now proceed to verify the invertibility condition of \cite[Th.\ 5]{nonsmooth_implicit},
	where we employ the Clarke Jacobian \( \clarke f \) as a conservative
	Jacobian \( \conserv f \).
	If \( \cys(\cdot) \) is differentiable at \( x \), \( \conserv f(x, \cy^{\star}(x)) = \left\{ [- \jac_1 h(x, \cy^{\star}(x)),~ I - \jac_2 h(x, \cy^{\star}(x)) ]\right\} \) and
	otherwise we have \( \conserv f(x, \cy^{\star}(x)) = \conv \{ [ - R_{i,1}, ~ I - R_{i,2}] \, | \, i \in \mathcal{P}(x) \} \).
	The mapping \( h(x, \cdot) \) is contractive which implies
	that \( \norm{\jac_2 h(x, \cy)} < 1 \) and \( I - \jac_2 h(x, \cy)  \) is invertible
	for any \( x \in \mathcal{X}, \cy \in \mathcal{Y}(x) \).
	Moreover, any convex combination of the contractive 
	linear maps \( R_{i,2} \), i.e.,
	\begin{equation*}
		\sum_{i \in \mathcal{P}(x)} \theta_i R_{i,2}, ~ \text{ with} ~ \theta_i \geq 0, ~ \text{and} ~ \sum_{i \in \mathcal{P}(x)} \theta_i =1
	\end{equation*}
	is also contractive.
	Thus, for any \( [A,\, B] \in \conserv f(x) \) the matrix \( B \) is invertible
	and the preconditions of \cite[Th.\ 5]{nonsmooth_implicit} are satisfied.
	We deduce that, for any \( x \in \mathcal{X} \), the conservative Jacobian of \( \cy^{\star}(\cdot) \) is
	given by
	\begin{equation} \label{eq:nonsmooth_sensitivity}
		\conserv \cys : x \rightrightarrows \{ - B^{-1} A : [A, ~ B] \in \conserv f (x, \cy^{\star}(x)) \}. 
	\end{equation}
	
	Next, we let \( S_1 = (S_{1,i})_{i \in \agents}, S_2 = (S_{2,i})_{i \in \agents} \) in 
	\autoref{alg:nonsmooth_sensitivity_learning_lqgs}, and notice that when 
	\( \norm{\tilde{\cy}^{\ell+1} - \tilde{\cy}^{\ell}} < \sigma \),
	this implies that \( S_1 = R_{j,1}, S_2 = R_{j,2} \) for some \( j \in \ml{P}(x) \).
	Therefore, the iteration \( \tilde{s}^{\ell+1} = S_2 \tilde{s}^{\ell} + S_1 \) is a contraction and
	converges to 
	\begin{equation} \label{eq:solution_map_jac_selection}
		\psi := (I - S_2)^{-1} S_1.
	\end{equation}
	We conclude by noting that \( \psi \in \conserv \cys(x) \)
	by \eqref{eq:nonsmooth_sensitivity}.
	\hfill \( \blacksquare \)

	\subsubsection{Proof of \autoref{lemma:s_bound_lqg}}
	\label{proof:s_bound_lqg}
	For sufficiently small \( \sigma \) and 
	\( \com(x, \tilde{\cy}^\ell) \in \tilde{\ml{Y}}_i, \, i \in \ml{P}(x) \),
	it holds true that \( [S_1, \, S_2] = R_i \).
	Therefore, the iteration \( \tilde{s}^{\ell+1} = S_2 \tilde{s}^\ell + S_1 \)
	is a contraction that converges to \( \psi = (I - S_2)^{-1} S_1 \in \conserv \cys(x) \), 
        see \autoref{prop:s_lqg_convergence}.
	Then, the following bounds hold:
	\begin{align*}
		\dist(\tilde{s}^\ell; \, \conserv \cys(x)) & = \min_{\chi \in \conserv \cys(x)} \norm{\tilde{s}^\ell - \chi} \\
		& \leq \norm{\tilde{s}^\ell - \psi} \\
		& \leq \frac{1}{1 - \eta} \norm{\tilde{s}^{\ell+1} - \tilde{s}^\ell},
	\end{align*}
	where the last inequality is by \cite[Th.\ 1.50(v)]{bauschke2017}.
	\hfill \( \blacksquare \)

	\subsubsection{Proof of \autoref{prop:s_lqgs_convergence}}
	\label{proof:s_lqgs_convergence} \( \jac \cys(\cdot) \) is constant since \( \cys(\cdot) \) is an affine mapping.
	The sequence \( \{\tilde{s}^\ell\}_{\ell \in \setN} \) converges linearly to \( \jac \cys(x) \), 
	for any \( x \in \ml{X} \), since \eqref{eq:s_update_nominal} and \eqref{eq:s_update_perturbed}
	coincide and the former is a contraction with \( \jac \cys(x) \) as its unique fixed point,
	as shown in \autoref{prop:s_nominal_convergence}.
	\hfill \( \blacksquare \)
	
		\subsubsection{Proof of \autoref{lemma:y_s_lqgs_bound}}
	\label{proof:y_s_lqgs_bound} 
	For the first bound we have
	\begin{align}
	& \norm{\cy^{k} - \cys(x^k)}  \overset{(a)}{=} \norm{h(x^{k}, \cy^{k-1}) - h(x^{k}, \cys(x^k))} \notag \\
		& \overset{(b)}{\leq} \eta \norm{\cy^{k-1} - \cys(x^k)} \notag \\
		& \overset{(c)}{\leq} \eta \norm{\cy^{k-1} - \cys(x^{k-1})} + \eta \norm{\cys(x^{k-1}) - \cys(x^k)} \notag \\
		& \overset{(d)}{\leq} \eta \norm{\cy^{k-1} - \cys(x^{k-1})} + \eta \norm{W} \norm{x^{k} - x^{k-1}}
		 \label{eq:y_bound_lqgs_sequential_intermediate} \\
		& \overset{(e)}{\leq} \eta^{k} \norm{\cy^0 - \cys(x^{0})} 
		+ \norm{W} \sum_{\ell=0}^{k-1} \eta^{k-\ell} \norm{x^{\ell+1} - x^{\ell}},
	\label{eq:y_bound_lqgs_sequential}
	\end{align}
	where (a) is by definition of \( \cy^{k} \), and \( \cys(x^{k}) \) being
	the unique fixed point of \( h(x^{k}, \cdot) \)
	(b) by contractiveness of \( h(x^{k}, \cdot) \) and
	(c) by adding/subtracting \( \cys(x^{k-1}) \) and the triangle inequality,
	(d) since \( \cys(\cdot) \) is affine, and
	(e) by telescoping.
	Further, we have that 
	\( \norm{x^{k+1} - x^{k}} = \beta^k \norm{\proj{\mathcal{X}}[x^{k} - \alpha^k \widehat{\nabla} \varphi_e^k] - x^k} \)
	and \( \proj{\mathcal{X}}[x^{k} - \alpha^k \widehat{\nabla} \varphi_e^k], x^k \in \mathcal{X} \),
	for all \( k \in \setN \), by convexity of \( \ml{X} \).
	Moreover, compactness of \( \ml{X} \) implies that there exists a constant \( B_{\mathcal{X}} > 0 \)
	such that \( \norm{x_1 - x_2} \leq B_{\mathcal{X}} \), for any \( x_1, x_2 \in \ml{X} \),
	which implies that \( \norm{x^{k+1} - x^{k}} \leq \beta^k B_{\ml{X}} \).
	Hence, substituting the last inequality onto \eqref{eq:y_bound_lqgs_sequential} we obtain
	\begin{equation*}
		\norm{\cy^{k} - \cys(x^k)} \leq \eta^{k} \norm{\cy^0 - \cys(x^0)} 
		+ \norm{W} B_{\mathcal{X}} \sum_{\ell=0}^{k-1} \eta^{k-\ell} \beta^k.
	\end{equation*}
	Letting \( B_{Y1} := \norm{\cy^0 - \cys(x^0)} \) and \( B_{Y2} := \norm{W} B_{\mathcal{X}} \)
	shows the first bound.
	
	For the second bound, contractiveness of \eqref{eq:s_update_perturbed} implies
	that
	\begin{equation}
		\norm{s^k - W} \leq \eta^k \norm{s^0 - W},
	\end{equation}
    and letting \( B_S := \norm{s^0 - W} \) establishes the claim.
	\hfill \( \blacksquare \)

	\subsection{Proofs of Section~\ref{subsec:outer_loop_an}}

	\subsubsection{Proof of \autoref{lemma:leader_is_definable}}
	\label{proof:leader_is_definable}
	Definability is preserved by function composition and, therefore,
	to show that \( \varphi_e \) is definable it suffices to
	prove that \( \cys(\cdot) \) is definable, see \cite[Exercise 1.11]{coste2000introduction}.
	To this end, we invoke \cite[Th.\ 5]{nonsmooth_implicit}
	since \( \cys(\cdot) \) is implicitly determined by
	\( f(x, \cy) = \cy - h(x, \cy) = 0 \).	
	Notice that, \( F \) is locally Lipschitz and definable
	as a consequence of \autoref{ass:followers_pg_general},
	while \( \proj{\ml{Y}(\cdot)}[\cdot] \) is locally Lipschitz and definable
	as it is PWA.
	Consequently, \( f \) is locally Lipschitz and definable and
	the invertibility condition of \cite[Th.\ 5]{nonsmooth_implicit}
	has been established in Appendix \ref{proof:s_lqg_convergence}.
	The claim follows by virtue of \cite[Th.\ 5]{nonsmooth_implicit} and
	uniqueness of the solution mapping.
	\hfill \( \blacksquare \)

	\subsubsection{Proof of \autoref{lemma:summable_errors_general}}
	\label{proof:summable_errors_general}
        We start by an auxiliary lemma that bounds \( \norm{e^k} \) w.r.t. the equilibrium and sensitivity error.
	\begin{lemma} \label{lemma:approximate_hypergradient_bound}
		Consider an iterate \( x^k \) of \eqref{eq:leaders_generic_iteration}
		and assume \( \ml{P}(x^k) \) is a singleton.
		  Then, there exist constants \( B_{ey}, B_{es} > 0 \) such that \( e^k \) satisfies the bound:
		\begin{equation} \label{eq:approximate_hypergradient_bound}
			\norm{e^k} \leq B_{ey} \norm{\cy^k - \cys(x^k)} + B_{es} \norm{s^k - \jac \cys(x^k)}.
		\end{equation}
	\end{lemma}
        \smallskip
        \begin{proof}
	Note that there exists \( B_{\cys}>0 \) such that \( \norm{\cys(x)} \leq B_{\cys} \),
	for any \( x \in \ml{X} \), and that \( \norm{\cy^k - \cys(x^k)} \) is upper-bounded by 
	\( \max_{k \in \setN} \sigma^k \).
	Thus, \( \{\cy^k\}_{k \in \setN} \) is contained in a compact set, and there exists \( B_{J \varphi} >0 \) such that \( \norm{\nabla_2 \varphi(x^k, y^k)} \leq B_{J \varphi} \),
	for any \( k \in \setN \), by continuity of \( \nabla_2 \varphi \).
	
	Now, to establish the desired bound we note that \( \varphi_e \) is differentiable at \( x^k \) and,
	therefore, we have
	\begin{align*}
		& \norm{e^k}  = \frac{1}{\alpha^k} \Big\lVert\proj{\ml{X}} \big[x^k - \alpha^k \widehat{\nabla} \varphi_e^k \big] - 
		\proj{\ml{X}} \big[x^k - \alpha^k \nabla \varphi_e(x^k) \big] \Big\rVert \\
		& \overset{(a)}{\leq} \norm{\widehat{\nabla} \varphi_e^k - \nabla \varphi_e(x^k)} \\
		& \overset{(b)}{=} \lVert \nabla_1 \varphi(x^k, \cy^k) + (s^k)^{\top} \nabla_2 \varphi(x^k, \cy^k) \\
		& ~~~ - \nabla_1 \varphi(x^k, \cys(x^k)) - (\jac \cys(x^k))^{\top} \nabla_2 \varphi(x^k, \cys(x^k)) \rVert \\
		& \overset{(c)}{\leq} \norm{\nabla_1 \varphi(x^k, \cy^k) - \nabla_1 \varphi(x^k, \cys(x^k))} + \\
		& ~~~~ \lVert (s^k - \jac \cys(x^k) )^{\top} \nabla_2 \varphi(x^k, \cy^k)  \\
		& ~~~~ + (\jac \cys(x^k))^{\top} (\nabla_2 \varphi(x^k, \cy^k) - \nabla_2 \varphi(x^k, \cys(x^k))) \rVert \\
		& \overset{(d)}{\leq} L_{\varphi 1} \norm{\cy^k - \cys(x^k)} 
		+ \norm{\nabla_2 \varphi(x^k, \cy^k)} \norm{s^k - \jac \cys(x^k)} \\
		& + \norm{\nabla_2 \varphi(x^k, \cy^k) - \nabla_2 \varphi(x^k, \cys(x^k)} \norm{\jac \cys(x^k)} \\
		& \overset{(e)}{\leq} (L_{\varphi 1} + L_S L_{\varphi 2})\norm{\cy^k - \cys(x^k)} 
		+ B_{J \varphi} \norm{s^k - \jac \cys(x^k)}
	\end{align*}
	where (a) is by nonexpansiveness of the projection,
	(b) by substituting the approximate and exact hypergradient, 
	(c) by adding/subtracting \( (\jac \cys(x^k))^{\top} \nabla_2 \varphi(x^k, \cy^k) \)
	and the triangle inequality,
	(d) by \autoref{ass:leaders_generic} and submultiplicativity of the matrix norm,
	(e) by \autoref{ass:leaders_generic} and the upper bounds on \( \norm{\nabla_2 \varphi(x^k, \cy^k)} \)
	and \( \norm{\jac \cys(x^k)} \).
        \end{proof}
        Next, we leverage the above result to prove \autoref{lemma:summable_errors_general}.
        
	\textit{Proof of \autoref{lemma:summable_errors_general}:} 
 By \autoref{ass:tolerance_distance}, we have that for any \( k > k' \) the bound \eqref{eq:s_bound_general} holds,
	thanks to the warmstarting of the inner loop.
	We apply \autoref{lemma:approximate_hypergradient_bound} by combining \eqref{eq:approximate_hypergradient_bound} with the error bounds
	\eqref{eq:y_bound} and \eqref{eq:s_bound_general} and the inner loop termination
	of \autoref{alg:innner_general_setup} (General) that yields
	\begin{align*}
		\alpha^k \norm{e^k} \leq \frac{B_{ey}}{1 - \eta} \alpha^k \sigma^k
		+ B_{es} B_{ps} \alpha^k \sigma^k + B_{es} B_{ns} \alpha^k \sigma^k.
	\end{align*}
	Invoking the design choice on \( \{\sigma^k\}_{k \in \setN} \)
	we deduce that \( \sum_{k=0}^{\infty} \alpha^k \norm{e^k} < \infty \).
	Notice that, if \( \ml{P}(x^k) \) is not a singleton for some \( k \in \setN \), then the corresponding sensitivity error  
        bound is not valid.
	Nonetheless, our summability result still holds since, by assumption, this happens for finitely many points.
	Finally, we complete the proof by recalling that any absolutely convergent sequence in a Banach space is 
	convergent \cite[Prob. 5.4.P6]{matrix_horn}, i.e., \( \sum_{k=0}^{\infty} \alpha^k \norm{e^k} < \infty \)
	implies that \( \sum_{k=0}^{\infty} \alpha^k e^k \to v \) for some vector \( v \in \setR^m \).
	\hfill \( \blacksquare \)
	
	\subsubsection{Proof of \autoref{th:general_convergence}}
	\label{proof:general_convergence}
	We let \( G(x) := - \conserv \varphi_e(x) - \ncone_{\ml{X}}(x)\) and observe that
        \( 0 \in G(x) \) implies that \( x \) is composite critical.
        \blue{To prove convergence, we will invoke \cite[Th.\ 3.2]{stochastic_sub_tame}. Hence, the remainder of the
        proof proceeds by showing that the assumptions in \cite{stochastic_sub_tame} hold true in our setting.}
	To make the paper self-contained, we recall these assumptions next.
	\begin{assumption}[\!\!{\cite[Assumption~A]{stochastic_sub_tame}}] \label{ass:general_inclusion_alg}
		\
		\begin{enumerate}
			\item All limit points of \( \{x^k\}_{k \in \setN}\) lie in \( \mathcal{X} \).
			\item The iterates are bounded, i.e.,
			\( \sup_{k \in \setN} \norm{x^k} < \infty \) and
			\( \sup_{k \in \setN} \norm{\xi^k} < \infty \).
			\item The sequence \( \{ \alpha^k\}_{k \in \setN}\) is nonnegative, nonsummable, and square-summable.
			\item The weighted noise sequence is convergent:
			\( \sum_{k=0}^{\infty} \alpha^k e^k \to v \) for some \( v \in \setR^m \).
			\item For any unbounded increasing sequence \( \{k_j\} \subseteq \setN \) such that
			\( x^{k_j} \) converges to some point \( \bar{x} \),
			\begin{equation*}
				\lim_{n \to \infty} \dist \bigg(\frac{1}{n} \sum_{j=1}^{n} \xi_{k_j};\, G(\bar{x}) \bigg) = 0.
				\tag*{\( \square \)}
			\end{equation*}
		\end{enumerate}
	\end{assumption}
	\begin{assumption}[\!\!{\cite[Assumption~B]{stochastic_sub_tame}}] \label{ass:lyapunov_condition}
		There exists a continuous function \( \phi : \setR^m \to \setR \),
		which is bounded from below, and such that the following two properties hold:
		\begin{enumerate}
			\item (Weak Sard) For a dense set of values \( r \in \setR \), the 
			intersection \( \phi^{-1}(r) \cap G^{-1}(0) \) is empty.
			\item (Descent) Whenever \( z : \setRp \to \setR^m \) is a trajectory of the differential
			inclusion \( \dot{z}(t) \in G(z(t)) \) and \( 0 \notin G(z(0)) \), there exists 
			a real \( T > 0 \) satisfying
			\begin{equation*}
				\phi(z(T)) < \sup_{t \in [0,T]} \phi(z(t)) \leq \phi(z(0)). \tag*{\( \square \)}
			\end{equation*}
		\end{enumerate}
	\end{assumption}
        
        \textit{Proof of \autoref{th:general_convergence}:}
        Observe that Assumption \ref{ass:general_inclusion_alg}.1
	and the \( x \)-part of \ref{ass:general_inclusion_alg}.2
	hold trivially by compactness of \( \mathcal{X} \).
	For the \( \xi \)-part of Assumption \ref{ass:general_inclusion_alg}.2
	we have that
	\begin{equation*}
		\norm{\xi^k} = \frac{1}{\alpha^k} \norm{\proj{\mathcal{X}}[x^k] - \proj{\mathcal{X}}[x^k - \alpha^k \zeta^k] }
		\leq \norm{\zeta^k} < \infty
	\end{equation*}
	where \( \zeta^k \in \conserv \varphi_e(x^k) \).
        Above, we used the fact that \( x^k \in \mathcal{X} \), for all \( k \in \setN \),
	which implies that \( x^k = \proj{\mathcal{X}}[x^k] \), 
	the nonexpansiveness of the projection and the boundedness of \( \nabla \varphi_e \).
	Further, Assumption \ref{ass:general_inclusion_alg}.3 holds by design, while
	Assumption~\ref{ass:general_inclusion_alg}.4 is shown in \autoref{lemma:summable_errors_general}.
	
	To show that \autoref{ass:general_inclusion_alg}.5 is satisfied, we first note that \( \conserv \varphi_e(\cdot) \) is convex-valued without loss
	of generality \cite{nonsmooth_implicit}, hence \( G(\cdot) \) is also convex-valued as the sum of
	two convex-valued mappings.
	Therefore, the following inequality holds:
	\begin{equation} \label{eq:ass_prox_distance_to_g}
		\dist \left(\frac{1}{n} \sum_{j=1}^{n} \xi_{k_j}, G(\bar{x})\right) \leq
		\frac{1}{n} \sum_{j=1}^{n} \dist \left( \xi_{k_j}, G(\bar{x})\right).
	\end{equation}
	Next, recall that the projected subgradient step can expressed as the solution
	to the following optimization problem
	\begin{equation*}
		\argmin_{w} \left\{\indic{\mathcal{X}}(w) + \frac{1}{2 \alpha^{k_j}} \norm{w - (x^{k_j}- \alpha^{k_j} \zeta^{k_j})}^2 \right\}.
	\end{equation*}
	Using Fermat's rule and the fact that \( \partial \indic{\mathcal{X}} = \ncone_{\mathcal{X}} \),
	we characterize the solution \( w^{k_j} \) of the previous problem as
	\begin{gather*}
		0  \in \ncone_{\mathcal{X}}(w^{k_j}) + \frac{1}{\alpha^{k_j}}(w^{k_j} - x^{k_j} + \alpha^{k_j} \zeta^{k_j}) \implies \\
		- \frac{1}{\alpha^{k_j}} (x^{k_j} - w^{k_j})  \in - \ncone_{\mathcal{X}}(w^{k_j}) - \zeta^{k_j}.
	\end{gather*}
	Note that \( w^{k_j} \to \bar{x} \) as \( j \to \infty \)
	since \( \alpha^{k_j} \to 0 \), and \( x^{k_j} \to \bar{x} \) by definition of \( \bar{x} \).
	Moreover, \( \conserv \varphi_e(\cdot) \) is outer semicontinuous, since it has 
	a closed graph and is locally bounded, and the same holds for \( \ncone_{\ml{X}} \).
	Then, we have that \( \xi^{k_j} = - \frac{1}{\alpha^{k_j}} (x^{k_j} - w^{k_j})  \) and,
	for some \( \bar{\zeta} \in \conserv \varphi_e(\bar{x}) \), it holds that
	\( \xi^{k_j} \to \ncone_{\mathcal{X}}(\bar{x}) + \bar{\zeta} \in G(\bar{x})\) .
	This implies that \( \dist(\xi^{k_j}, G(\bar{x})) \to 0 \) and, consequently,
	it holds \( \lim_{n \to \infty}\frac{1}{n} \sum_{j=1}^{n} \dist \left( \xi_{k_j}, G(\bar{x})\right) = 0 \).
	Using \eqref{eq:ass_prox_distance_to_g}, we establish \autoref{ass:general_inclusion_alg}.5.
	
	To verify \autoref{ass:lyapunov_condition}, we will utilize \( \varphi_e \) as the Lyapunov function
	\( \phi \) and recall that \( \varphi_e \) and \( \ml{X} \) are definable by virtue of 
	\autoref{lemma:leader_is_definable} and \autoref{ass:leaders_generic}, respectively.
	Therefore, \( \varphi_e \) and \( \ml{X} \) admit a chain rule as in \cite[Def.\ 5.1]{stochastic_sub_tame}
	according to \cite[Th.\ 5.8]{stochastic_sub_tame}.
	Then, \autoref{ass:lyapunov_condition}.1 follows exactly by the arguments in the proof
	of \cite[Cor.\ 6.4]{stochastic_sub_tame},
	while \autoref{ass:lyapunov_condition}.2 is verified by \cite[Lemma 6.3]{stochastic_sub_tame}.
	Having established its preconditions, the result follows by applying
	\cite[Th.\ 3.2]{stochastic_sub_tame}.
	\hfill \( \blacksquare \)

	\subsubsection{Proof of \autoref{th:convergence_lqg}}
	\label{proof:convergence_lqg}
        To show the claim, it suffices to prove summability of the errors in the presence
        of nonsmoothness.
	Then, the remainder of the proof is identical to the proof of \autoref{th:general_convergence}.
 	To this end, we will first extend the hypergradient error bound in \eqref{eq:approximate_hypergradient_bound}
	to the nonsmooth regime.
	According to the conservative Jacobian chain rule \cite[Prop.\ 2]{nonsmooth_implicit}
	the (conservative) hypergradient is given by
	\begin{align}
		& \conserv \varphi_e : x \rightrightarrows \{H_{\varphi_e, \chi} \, | \, \chi \in \conserv \cys(x) \}, \text{ where} \\
		& H_{\varphi_e, \chi}(x) := \nabla_1 \varphi (x, \cys(x)) + \chi^{\top} \nabla_2 \varphi(x, \cys(x)).
	\end{align}
	Then, the following bound holds
	\begin{align*}
		& \norm{e^k}  \overset{(a)}{=} \min_{\hat{\zeta}^k \in \conserv \varphi_e(x^k)} 
		\frac{1}{\alpha^k} \norm{\proj{\ml{X}}[x^k - \alpha^k \widehat{\nabla} \varphi_e^k] - 
		\proj{\ml{X}}[x^k - \alpha^k \hat{\zeta}^k]} \\
		& \overset{(b)}{\leq} \frac{1}{\alpha^k} 
		\norm{\proj{\ml{X}}[x^k - \alpha^k \widehat{\nabla} \varphi_e^k] -
		\proj{\ml{X}}[x^k - \alpha^k H_{\varphi_e, \psi}]} \\
		& \overset{(c)}{\leq} \norm{\widehat{\nabla} \varphi_e^k - H_{\varphi_e, \psi}(x^k) } \\
		& \overset{(d)}{=} \lVert \nabla_1 \varphi(x^k, \cy^k) + {s^k}^{\top} \nabla_2 \varphi(x^k, \cy^k) \\
		&	~~~ - \nabla_1 \varphi(x^k, \cys(x^k)) - \psi^{\top} \nabla_2 \varphi(x^k, \cys(x^k)) \rVert \\
		& \overset{(e)}{\leq} (L_{\varphi 1} + L_S L_{\varphi 2})\norm{\cy^k - \cys(x^k)} 
		+ B_{J \varphi} \norm{s^k - \psi},
	\end{align*}
	where (a) is by \eqref{eq:error_def_general},
	(b) by noting that \( \psi \in \conserv \cys(x^k) \) as in \eqref{eq:solution_map_jac_selection} 
	which implies that \( H_{\varphi_e, \psi} \in \conserv \varphi_e(x^k) \),
	(c) by nonexpansiveness of the projection,
	(d) by substitution,
	and (e) follows exactly as in the proof of \autoref{lemma:approximate_hypergradient_bound}.
	
	Combining the error bounds \eqref{eq:y_bound} and \eqref{eq:s_bound_lqg} and the termination
	criterion of \autoref{alg:nonsmooth_sensitivity_learning_lqgs}, 
        which together with \autoref{ass:tolerance_distance} guarantees that \eqref{eq:s_bound_lqg} is valid,
        we obtain
	\begin{equation}
		\alpha^k \norm{e^k} \leq \frac{L_{\varphi 1} + L_S L_{\varphi 2}}{1 - \eta} \alpha^k \sigma^k
		+ \frac{B_{J \varphi}}{1 - \eta} \alpha^k \sigma^k.
	\end{equation}
	Thus, the sequence \( \{a^k \norm{e^k}\}_{k \in \setN} \) is summable and the same holds
	for \( \{a^k e^k \}_{k \in \setN} \), which completes the proof.
	\hfill \( \blacksquare \)
	
	\subsubsection{Proof of \autoref{lemma:summable_errors_lqgs}}
	\label{proof:summable_errors_lqgs}
	It can be shown, as in \autoref{lemma:approximate_hypergradient_bound},
	that the normed error in \eqref{eq:error_def_lqgs} satisfies the bound
	\begin{equation} \label{eq:approximate_hypergradient_bound_lqgs}
		\norm{e^k} \leq \alpha B_{ey} \norm{\cy^k - \cys(x^k)} + \alpha B_{es} \norm{s^k - W}.
	\end{equation}
	Combining \eqref{eq:approximate_hypergradient_bound_lqgs} and \autoref{lemma:y_s_lqgs_bound}, we
	can derive the following
	\begin{equation}
		\beta^k \norm{e^k} \leq \delta_1 \beta^k \eta^k 
		+ \delta_2 \beta^k \sum_{\ell=0}^{k-1} \eta^{k-\ell} \beta^\ell
	\end{equation}
	for constants
	\( \delta_1 := (B_{ey} B_{Y1} + B_{es} B_S) \alpha \) and
	\( \delta_2 := B_{ey} B_{Y2} \alpha \).
	We will prove summability by considering each of the two terms individually.
	
	For the first term, we have 
	\( \sum_{k = 0}^{\infty} \beta^k \eta^k \leq \sum_{k = 0}^{\infty} \eta^k < \infty \),
	where the first inequality is since \( \beta^k \leq 1 \), while the second by \( \eta < 1 \).
	
	For the second term, it holds true that
	\begin{equation*}
		\sum_{k = 0}^{\infty} \beta^k \sum_{\ell=0}^{k-1} \eta^{k-\ell} \beta^\ell \leq 
		\sum_{k = 0}^{\infty} \sum_{\ell=0}^{k-1} \eta^{k-\ell} (\beta^\ell)^2 < \infty,
	\end{equation*}
	where the first inequality is due to \( \beta^k \leq \beta^\ell \) for all \( k \geq \ell \),
	while the second is by virtue of \( \sum_{k=0}^{\infty} (\beta^k)^2 < \infty \)
	and \cite[Lemma 3.1]{sequences_lemma},
	hence completing the proof.
	\hfill \( \blacksquare \)
	
	\subsubsection{Proof of \autoref{th:lqgs_convergence}}
	\label{proof:lqgs_convergence}
	The outer iterations of \autoref{alg:general_hypergradient}, with inner loop as in \autoref{alg:innner_general_setup} (LQSG),
	can be viewed as the inexact
	KM process \( x^{k+1} = x^k + \beta^k (T(x^k) + e^k - x^k) \),
	where \( T : x \mapsto \proj{\ml{X}}[x - \alpha^k \nabla \varphi_e(x)] \) and
	the errors \( e^k \) are as in \eqref{eq:error_def_lqgs}.
	For any \( \{ \alpha^k \}_{k \in \setN} \) as in \autoref{lemma:summable_errors_lqgs},
	the operator \( T(\cdot) \) is nonexpansive by virtue of 
         Theorem 18.15(v), Proposition 26.1(iv)(d), and Remark 4.34(i) in \cite{bauschke2017}.
	Moreover, the set \( \fix(T) \) coincides with the set of solutions to \eqref{eq:upper_level_substi},
	which is nonempty by continuity of \( \varphi_e \) and compactness of \( \ml{X} \).
	
	Next, we note that \( \sum_{k = 0}^{\infty} \beta^k \norm{e^k} < \infty \) by \autoref{lemma:summable_errors_lqgs}
	and \( \sum_{k = 0}^{\infty} \beta^k(1 - \beta^k) = \infty \).
	Hence, the preconditions of \cite[Prop.\ 5.34]{bauschke2017} are met,
	which implies that \( \{x^k\}_{k \in \setN} \) converges to \( \fix(T) \)
	and the claim follows.
	\hfill \( \blacksquare \)
	
	\subsubsection{Proof of \autoref{th:lqgs_convergence_strong}}
	\label{proof:lqgs_convergence_strong}
	Initially, note that \( \varphi_e(x) = \varphi([I;\, W] x + [0;\, w]) \) and
	the matrix \( [I;\, W] \) is full column rank, hence, \( \varphi_e \)
	is \( \sigma_{\varphi_e} \)-strongly convex,
	with \( \sigma_{\varphi_e}:= \sigma_{\varphi} \sigma_{\min}^2 \)
	and \( \sigma_{\min} > 0 \) being the smallest singular value of \( [I;\, W] \).
	This implies that \eqref{eq:upper_level_substi} admits a unique solution, 
	denoted \( x^{\star} \).
	Further, we recall that projected gradient descent converges linearly for	strongly convex problems \cite[Prop.~26.16]{bauschke2017}, i.e., 
	\begin{equation} \label{eq:linear_convergence_strongly}
		\norm{\tilde{x}^k - x^{\star}} \leq (1 - \kappa(\alpha)) \norm{x^k - x^{\star}}
	\end{equation}
	where \( \tilde{x}^k := \proj{\ml{X}}[x^k - \alpha \nabla \varphi_e(x^k)] \) is the
	exact projected gradient step,
	and \( \kappa(\alpha) := 1 - \sqrt{1 - \alpha \sigma_{\varphi_e}} \).
	We exploit this property to derive useful error bounds on
	\( \norm{x^k - x^{\star}} \) and \( \norm{\cy^k - \cys(x^k)} \).

	For brevity, we let \( \hat{x}^k := \proj{\ml{X}}[x^k - \alpha \widehat{\nabla} \varphi_e^k] \)
	denote the inexact projected gradient step, and for the leader's error we have
	\begin{align}
		& \norm{x^{k+1} - x^{\star}} 
		 \overset{(a)}{=} \norm{x^k + \beta^k(\hat{x}^k - x^k) - x^{\star}} \notag \\
		& \overset{(b)}{\leq} (1 - \beta^k) \norm{x^k - x^{\star}} + \beta^k \norm{\hat{x}^k - x^{\star}} 
		\notag \\
		& \overset{(c)}{\leq} (1 - \beta^k) \norm{x^k - x^{\star}} + \beta^k \norm{\hat{x}^k - \tilde{x}^k} 
		+ \beta^k \norm{\tilde{x}^k - x^{\star}} \notag \\
		& \overset{(d)}{\leq} \left(1 - \beta^k \kappa(\alpha) \right) \norm{x^k - x^{\star}} 
		+ \beta^k \norm{\hat{x}^k - \tilde{x}^k} \notag \\		
		& \overset{(e)}{\leq} \left(1 - \beta^k \kappa(\alpha) \right) \norm{x^k - x^{\star}} \notag \\
		& ~~~~ + \beta^k \alpha B_{ey} \norm{\cy^k - \cys(x^k)} + \beta^k \alpha B_{es} \norm{s^k - W}
		\label{eq:x_bound_lqgs}
	\end{align}
	where (a) is by substituting the update rule for \( x^{k+1} \),
	(b) by adding/subtracting \( \beta^k x^{\star} \) and the triangle inequality,
	(c) by the adding/subtracting \( \tilde{x}^k \) and triangle inequality,
	(d) by \eqref{eq:linear_convergence_strongly},
	(e) by observing that \( e^k =\hat{x}^k - \tilde{x}^k \) and using \eqref{eq:approximate_hypergradient_bound_lqgs}.
	
	For the NE error we have
	\begin{align}
		& \norm{\cy^{k+1} - \cys(x^{k+1})}  \overset{(a)}{\leq} \eta \norm{\cy^k - \cys(x^{k})} \notag \\
		& ~~~~~ + \eta \norm{W} \norm{x^k - x^{k+1}} \notag \\
		& \overset{(b)}{\leq} \eta \norm{\cy^k - \cys(x^{k})} + \beta^k \eta \norm{W} \norm{\hat{x}^k - x^{\star}} \notag \\
		& ~~~~~  + \beta^k \eta \norm{W} \norm{x^k - x^{\star}} \notag \\
		& \overset{(c)}{\leq} 
		\left(\eta + \beta^k \eta \norm{W} \alpha B_{ey} \right) \norm{\cy^k - \cys(x^{k})} \notag \\
		& ~~~~~ + \beta^k \eta \norm{W} B_{es} \alpha \norm{s^k - W} \notag \\
		& ~~~~~ + \beta^k \eta \norm{W} \left(2 - \kappa(\alpha) \right) \norm{x^k - x^{\star}}
		\label{eq:y_bound_lqgs}
	\end{align}
	where (a) is shown in \eqref{eq:y_bound_lqgs_sequential_intermediate}, 
	(b) by substituting the update rule of \( x^{k+1} \), 
	adding/subtracting \( x^{\star} \) and the triangle inequality,
	and (c) by bounding \( \norm{\hat{x}^k - x^{\star}} \) as done in \eqref{eq:x_bound_lqgs}.
	
	Next, we let \( \varepsilon^{k} := \left(\norm{x^{k} - x^{\star}}, \norm{y^{k} - \cys(x^{k})}, \norm{s^{k} - W}\right) \)
	and compactly express the bounds \eqref{eq:sens_vanish_bound}, \eqref{eq:x_bound_lqgs}, and \eqref{eq:y_bound_lqgs} as
	\begin{equation*}
		\norm{\varepsilon^{k+1}} \leq M(\beta^k) \norm{\varepsilon^k},
	\end{equation*}
	where
	\begin{align}
		& M(\beta^k) := 
		\diag(1, \eta, \eta) + \beta^k E \notag \\
		& E :=
		\begin{bmatrix}
				- \kappa(\alpha) & \alpha B_{ey} & \alpha B_{es} \\
				(2 - \kappa(\alpha)) \eta \norm{W} &  
				\eta \alpha \norm{W} B_{ey} & \eta \alpha \norm{W} B_{es}  \notag \\
				0 & 0 & 0
		\end{bmatrix}
		.
	\end{align}
	We will establish that \( \varepsilon^k \) vanishes at a linear rate by proving that 
	\( M(\beta) \) has a spectral radius \( \rho(M(\beta)) < 1 \) for any sufficiently small \( \beta > 0 \).
        This implies that the sequence \( \{ (x^k, \cy^k, s^k) \}_{k \in \setN} \) converges to \( (x^{\star}, \cys(x^{\star}), \jac \cys(x^{\star})) \),
        and thus it proves the claim.
	
	In that direction, we note that, for any small \( \beta > 0 \) the matrix \( M(\beta) \) is non-negative and irreducible \cite[Def.\ 1.3.16]{matrix_horn},
	since the associated directed graph is strongly connected \cite[Th.\ 6.2.24]{matrix_horn}.
	Further, its diagonal is nonzero implying that it is a primitive matrix \cite[Def.\ 8.5.0]{matrix_horn} by virtue
	of \cite[Th.\ 8.5.2]{matrix_horn} and \cite[Th.\ 8.5.9]{matrix_horn}.
	Therefore, the Perron-Frobenius theorem \cite[Th.\ 8.4.4]{matrix_horn} ensures that the spectral radius \( \rho(M(\beta)) \) is
	equal to the maximum real eigenvalue \( \lambda_{\max} \) of \( M(\beta) \),
	and all other eigenvalues have strictly smaller modulus.
	
	Observe that \( \lambda_1 = 1 \) is a simple eigenvalue of \( M(0) \) and \( v = (1,0,0) \) is both the left and right corresponding eigenvector,
	while \( \lambda_2 = \eta < 1 \) is the other eigenvalue, of multiplicity 2;
	hence, \( \rho(M(0)) = \lambda_1 \).
	Then, invoking \cite[Th.\ 1]{eigen_perturbation} we have that the derivative of \( \lambda_1 \) with respect to \( \beta \) satisfies
	\( \nabla_{\beta} \lambda_1 = v^{\top} E v = - \kappa(\alpha) < 0 \).
	Therefore, by continuity of the eigenvalues with respect to \( \beta \) we deduce that for \( \beta > 0 \)
	sufficiently small it holds \( \rho(M(\beta)) < 1 \).
	Further, by direct computation we have that \( 1  \) is an eigenvalue of \( M(\beta) \),  if and only if
	\( \beta \) is equal to zero or \( \beta \) is equal to
	\begin{equation}
		\overline{\beta} := \frac{\kappa(\alpha)(1 - \eta)}
		{2 B_{ey} \alpha \eta \norm{W}}.
	\end{equation}
	Finally, we deduce that for any step size sequence \( \{\beta^k\}_{k \in \setN} \) that is identically equal to \( \beta \in (0, \overline{\beta}) \) the error
	sequence \( \{ \varepsilon^k \}_{k \in \setN} \) is driven to zero at a linear rate upper-bounded by \( \rho(M(\beta)) \), which concludes the proof.
	\hfill \( \blacksquare \)	
\balance
\bibliographystyle{ieeetr}
\bibliography{references}

\end{document}